%% file: main.tex
\documentclass[twocolumn,letterpaper]{IEEEAerospaceCLS}  

\usepackage[]{graphicx}    
\newcommand{\ignore}[1]{}  

\usepackage{graphbox}
\usepackage{array}
\usepackage{amssymb}
\usepackage[version=4]{mhchem}
\usepackage{siunitx}
\usepackage{longtable,tabularx}
\usepackage{booktabs}
\usepackage{mathtools}
\usepackage{multirow}
\usepackage{caption}
\usepackage{subcaption}
\usepackage{makecell}
\usepackage{booktabs}

\usepackage{amsthm}

\usepackage[T1]{fontenc}\usepackage{tikz}
\usepackage{amsmath, amsfonts}
\usetikzlibrary{arrows.meta}

\usepackage[linesnumbered,ruled,vlined,norelsize,noend]{algorithm2e}
\SetInd{5pt}{10pt}
\SetNlSty{relax}{}{\enspace}

\def\thickhline{\noalign{\hrule height.8pt}}

\usepackage{xcolor}
\definecolor{light-gray}{gray}{0.95}

\usepackage{hyperref}

\usepackage[font=small,skip=4pt]{caption}
\begin{document}
\title{Agile Tradespace Exploration for Space Rendezvous Mission Design via Transformers}

\author{%
Yuji Takubo\\ 
Dept. of Aeronautics \& Astronautics\\
Stanford University\\
496 Lomita Mall, Stanford, CA 94305\\
ytakubo@stanford.edu
\and 
Daniele Gammelli\\
Dept. of Aeronautics \& Astronautics\\
Stanford University\\
496 Lomita Mall, Stanford, CA 94305\\
gammelli@stanford.edu
\and 
Marco Pavone \\ 
Dept. of Aeronautics \& Astronautics\\
Stanford University\\
496 Lomita Mall, Stanford, CA 94305\\
pavone@stanford.edu
\and
Simone D'Amico \\ 
Dept. of Aeronautics \& Astronautics\\
Stanford University\\
496 Lomita Mall, Stanford, CA 94305\\
damicos@stanford.edu
}

\maketitle

\thispagestyle{plain}
\pagestyle{plain}

\maketitle

\thispagestyle{plain}
\pagestyle{plain}

\setcounter{footnote}{0}

\begin{abstract}
Spacecraft rendezvous enables on-orbit servicing, debris removal, and crewed docking, forming the foundation for a scalable space economy. 
Designing such missions requires rapid exploration of the tradespace between control cost and flight time across multiple candidate targets.
However, multi-objective optimization in this setting is challenging, as the underlying constraints are often nonconvex, and mission designers must balance accuracy (e.g., solving the full problem) with efficiency (e.g., convex relaxations), slowing iteration and limiting design agility.
To address these challenges, this paper proposes an AI-powered framework that enables agile and generalized rendezvous mission design. 
Given the orbital information of the target spacecraft, boundary conditions of the servicer, and a range of flight times, a transformer model generates a set of near-Pareto optimal trajectories across varying flight times in a single parallelized inference step, thereby enabling rapid mission trade studies. 
The model is further extended to accommodate variable flight times and perturbed orbital dynamics, supporting realistic multi-objective trade-offs.
Validation on chance-constrained rendezvous problems in Earth orbits with passive safety constraints demonstrates that the model generalizes across both flight times and dynamics, consistently providing high-quality initial guesses that converge to superior solutions in fewer iterations.
Moreover, the framework efficiently approximates the Pareto front, achieving runtimes comparable to convex relaxation by exploiting parallelized inference.
Together, these results position the proposed framework as a practical surrogate for nonconvex trajectory generation and mark an important step toward AI-driven trajectory design for accelerating preliminary mission planning in real-world rendezvous applications.
\end{abstract} 

\tableofcontents

\section{Introduction}

\vspace{1.1em}

\begin{figure}[h!]
    \centering
    \includegraphics[width=\columnwidth]{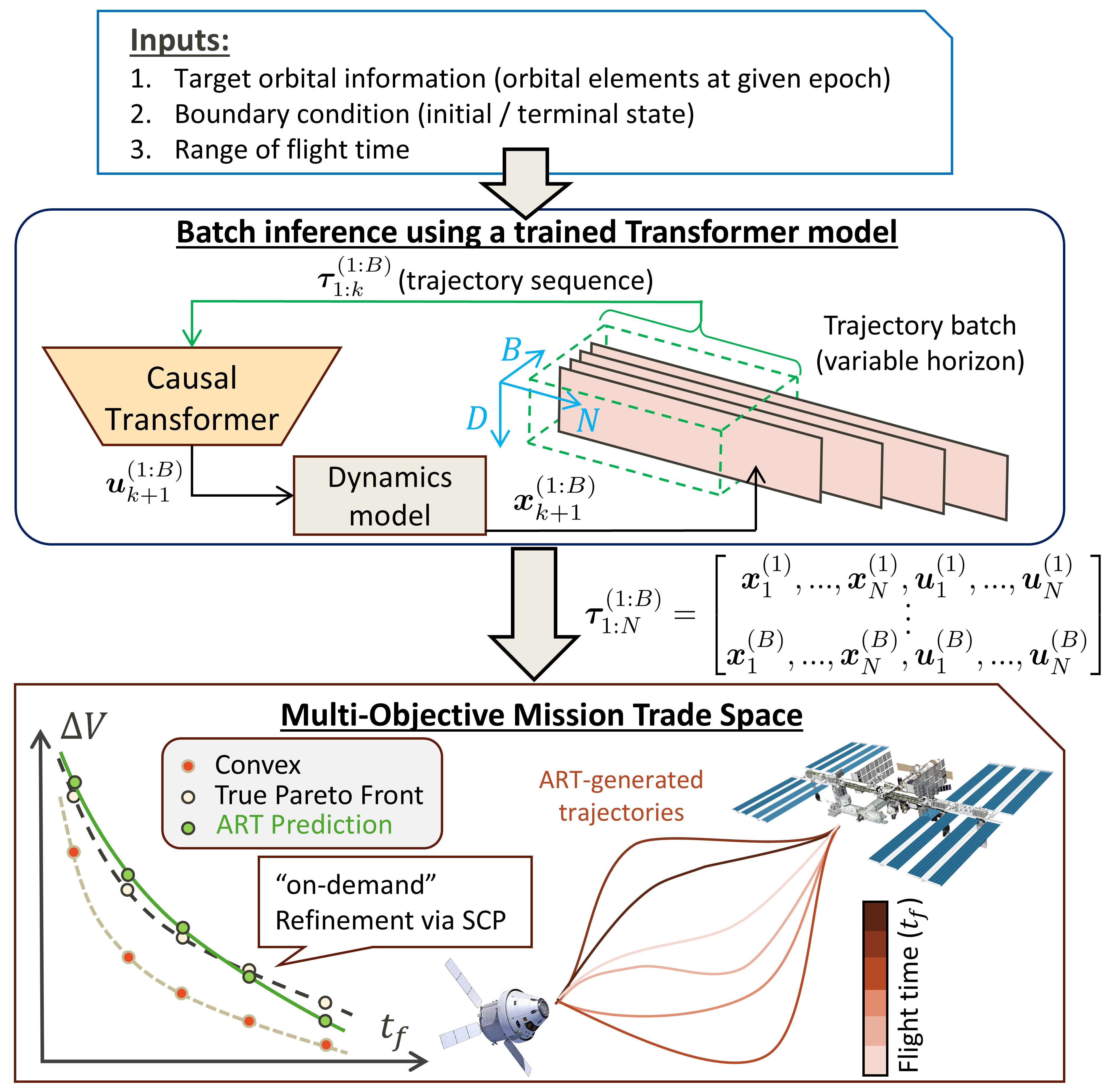}
    \caption{This paper presents a framework for AI-powered agile mission design. Mission configurations (target orbits, boundary conditions, and flight-time ranges) serve as inputs to a trajectory generation pipeline (top), where a transformer model autoregressively predicts control sequences that are propagated through the system dynamics (middle). Batch inference efficiently produces multiple Pareto-optimal trajectories across flight times that are encoded into a tensor of $[B (\text{batch size}),N(\text{time horizon}),D (\text{types of subsequence elements})]$, which can be projected onto the multi-objective solution space as surrogate mission models. Selected point solutions can then be further refined using sequential convex programming (bottom). }
    \label{fig:art_md_summary}
\end{figure}

Distributed Space Systems (DSS) offer capabilities that extend far beyond those of traditional monolithic spacecraft, enabling a new class of missions and operational concepts.
In recent years, DSS architectures have supported a wide range of applications, including fundamental science missions \cite{montenbruck2006vector}\cite{montenbruck2008navigation}\cite{Guffanti2023VISORS}\cite{lowe2023reduced}, crewed rendezvous operations \cite{dsouza2007orion}, and potential enhancements in space situational awareness \cite{kruger2025orbit}\cite{Rizza2025FALCON}. 
A particularly transformative development lies in high-precision rendezvous and docking technologies, which are redefining the paradigm of space logistics and on-orbit servicing. 
These advances enable critical functions such as refueling, repair, and debris removal, thereby extending spacecraft lifetimes and enhancing the sustainability of space operations. 
Looking ahead, logistics vehicles will be expected to perform frequent Rendezvous, Proximity Operations, and Docking (RPOD) maneuvers with multiple client spacecraft or objects dispersed across various orbital domains around Earth and cislunar space. 
Each potential servicing target introduces a unique set of orbital characteristics and phasing opportunities.
Consequently, mission designers must navigate a complex, combinatorial design space---rapidly iterating through catalogs of candidate mission profiles to identify feasible trajectories and select the most advantageous target for a given reference orbit---underscoring both the opportunities and challenges inherent in deploying DSS for next-generation space logistics.

Historically, space mission design has been an arduous, highly manual process, largely dependent on expert knowledge and rarely regarded as a component of onboard autonomy. 
This difficulty arises from two main factors: (i) its formulation as a multi-objective optimization problem, in which several competing objectives (e.g., flight time, fuel consumption, radiation exposure, observability) must be optimized simultaneously, and (ii) the presence of complex nonconvex and often implicit constraints, such as safety and observation requirements. 
The result is a constrained Multi-Objective Optimal Control Problem (MO-OCP), where the Pareto front of feasible solutions is of central importance. 
Because general MO-OCPs are extremely challenging to solve, traditional practice has relied on problem-specific search strategies \cite{landau2022star} or metaheuristics \cite{li2010optimal}, typically under strong simplifying assumptions, to produce preliminary trajectories of acceptable quality.
However, while useful in early mission design, these methods suffer from limited scalability and generalizability. Moreover, they require a full recomputation whenever mission parameters change, underscoring their fragility in dynamic design environments.

To address these challenges, recent work has investigated learning-based methods that can efficiently generate spacecraft trajectories from a single trained model. 
Learning-based techniques have been applied to diverse scenarios, from interplanetary transfers to rendezvous. 
Yet, for a mission design tool suitable for operational use, such as that envisioned for RPOD missions, several critical capabilities remain underdeveloped. 
First, generalization is essential; the governing dynamics depend on the reference orbit \cite{damico_phd_2010}\cite{hill1878researches}\cite{yamanaka2002new}, and models must therefore be able to extrapolate trajectory generation to previously unseen relative dynamics.
Second, limited attention has been devoted to multi-objective offline learning frameworks \cite{hayes2021practical}\cite{zhu2023scaling}, which are indispensable for practical design applications that require trade-space exploration.

To address these limitations, this study introduces an agile transformer-based trajectory generation pipeline for RPOD mission design (cf. Figure \ref{fig:art_md_summary}).
This framework is capable of generating batches of near-Pareto-optimal and near-feasible trajectories in a single inference pass.
This approach builds on the Autonomous Rendezvous Transformer (ART) framework \cite{guffanti2024transformers}\cite{celestini2024transformer}\cite{celestini2025generalizable}\cite{takubo2025towards}\cite{takubo2026semantic}, originally developed as a warm-starting technique for nonconvex trajectory optimization to accelerate onboard guidance. 
In this work, ART is extended and repurposed for mission design, specifically for conducting trade studies in complex scenarios and for approximating the Pareto front of nonconvex RPOD trajectory optimization problems.

The contributions of this paper are threefold: 
\begin{itemize}
    \item This paper proposes a framework for AI-powered (transformer-based) agile space mission design that generates a near-Pareto optimal solution set for a given mission configuration.
    \item This paper extends the ART architecture to meet the requirements of space mission design. Specifically, this work investigates (i) the generalization across dynamics and flight time, and (ii) the generation of Pareto-optimal trajectories via parallelizable inference. 
    \item Numerical experiments validate these requirements through an extensive set of rendezvous mission profiles in perturbed-Keplerian orbits, incorporating data from the latest space debris catalog. Crucially, trajectory generation with ART achieves runtimes comparable to those of convex relaxation methods, underscoring its effectiveness as a practical and agile tool for multi-objective nonconvex trajectory generation.
\end{itemize}

The remainder of this paper is organized as follows. 
Section~\ref{sec:related_works} reviews prior work on mission design and learning-based trajectory generation, while Section~\ref{sec:background} introduces the necessary background on ART.
Section \ref{sec:methodology} presents the formulation of the MO-OCP and introduces the ART extension for Pareto-optimal trajectory generation. 
Section \ref{sec:case_study} formulates the spacecraft rendezvous problem in perturbed Keplerian orbits as a chance-constrained nonconvex optimization problem. 
Section \ref{sec:results} presents the generalization performance of the proposed framework and analyzes the resulting multi-objective solution surface.
Finally, Section \ref{sec:conclusion} concludes the paper and discusses directions for future work.

\section{Related Work} \label{sec:related_works}

This work addresses the challenges of preliminary trajectory design in space missions by bridging traditional methodologies with emerging learning-based trajectory generation approaches.

\subsection{Space Mission Design}

Space mission design remains a manual and computationally intensive process due to the complexity of multi-objective trades under stringent constraints. Existing methods can broadly be classified into three categories.

\subsubsection{1. Search-based heuristics} One of the popular approaches is the deterministic broad search under problem-specific simplifying assumptions (e.g., patched-conics with impulsive maneuvers \cite{landau2022star}\cite{takubo2024automated}, tree-based motion planning techniques \cite{starek2017fast}\cite{deka2023astrodynamics}). 
While state-of-the-art algorithms enable a fast broad search, crude approximation leads to a non-negligible gap between the approximated objective function during the initial guess search and the true objective function in the later refinement \cite{takubo2026preliminary}. 
Another common approach is using stochastic, population-based methods such as genetic algorithms and differential evolution \cite{li2010optimal}\cite{izzo2010global}\cite{vasile2011multi}\cite{englander2017automated}\cite{jin2020robust}. 
While capable of exploring global design spaces and accommodating uncertainty in high-fidelity simulations \cite{Arya2024Stochastic}\cite{takubo2022robust}\cite{geller2022robust}, these methods generally require extensive tuning, suffer from poor scalability, and provide limited generalization.

\subsubsection{2. Scalarization of objectives} Another strategy is to reduce the multi-objective problem to a scalarized one, either through weighted sums or more sophisticated scalarization schemes \cite{vasile2019multi}.
Such formulations can then be solved via nonconvex optimization methods \cite{berning2024chance}\cite{margolis2024robust}, often leveraging algorithms such as Sequential Convex Programming (SCP) \cite{malyuta_scp_2022}. 
Hybrid methods combining stochastic exploration with gradient-based local refinement have also been explored \cite{vasile2019multi}.
However, approximating the Pareto set typically requires sweeping large sets of weights or thresholds, resulting in high computational cost.
Moreover, incorporating expressive constraints, such as signal temporal logic or mixed-integer constraints \cite{malyuta2023fast}\cite{szmuk2020successive}\cite{mao2022successive}, further increases the risk of non-convergence.

\subsubsection{3. Reduced-order and surrogate models} Simplified dynamics and surrogate models are often employed to accelerate the design of DSS \cite{lowe2023reduced}\cite{lippe2021phd}.
While these methods reduce computational burden, they may neglect operational constraints (e.g., passive safety \cite{damico_phd_2010}\cite{aguilar2022abort}\cite{guffanti_jgcd_2023}\cite{takubo2025safe}, field-of-view/approach cones \cite{kim2010convex}\cite{lee2014dual}\cite{takubo2024multiplicative}, plume impingement effects \cite{malyuta2023fast}, etc.) that strongly influence feasibility. 
As a result, preliminary designs generated under these simplifications often require large corrections. 
Crucially, these methods typically do not yield full trajectories but only approximate mission profiles. 

In summary, there remains a gap between computationally intensive, high-fidelity trajectory design methods and the fast but coarse approximations provided by surrogate models. 
This paper seeks to bridge that gap by introducing an AI–based trajectory generation framework that accelerates the design process while preserving generality.  
The proposed approach is highly flexible and capable of predicting accurate Pareto solution surfaces under complex constraints, without requiring re-training the neural network model when mission configurations change.

\subsection{Learning-based Trajectory Generation in Space}

In recent years, there has been a growing interest in applying machine learning and AI techniques to generate spacecraft trajectories for both mission design and onboard autonomy contexts. 
This is primarily motivated by the goal of facilitating trajectory generation by (i) reducing the overall computational cost, and (ii) enabling global search over candidate trajectories, which conventional methods often struggle to achieve without extensive expert hand-tuning.
Learning-based methods have been applied in different capacities, including the prediction of optimal co-states for indirect methods \cite{izzo2021real}, identifying active constraint sets \cite{briden2025transformer}, or estimating primal-dual variables \cite{briden2025constraint}. 
Direct trajectory generation has also been pursued via imitation learning \cite{tsukamoto2024neural} and reinforcement learning \cite{broida2019spacecraft}\cite{zavoli2021reinforcement}, leveraging both transformer- \cite{guffanti2024transformers} and diffusion-based architectures \cite{briden2025compositional}\cite{li2024efficient}. 
Neural architectures have been applied to highly multi-modal trajectory optimization problems, such as Earth–Moon low-thrust transfers \cite{Li2023Amortized} \cite{graebner2025global}, including cases with missed-thrust contingencies \cite{sinha2025initial}.
Many of these methods are used to provide warm-starts \cite{Banerjee_2020}\cite{yuan2025filtering} for nonlinear optimization problems, which have recently been demonstrated on Astrobee aboard the International Space Station \cite{banerjee2025deep}.
Within this context, ART has extended transformer-based trajectory generation to a variety of settings, including model predictive control \cite{celestini2024transformer}, chance-constrained and fault-tolerant optimal control \cite{takubo2025towards}, and trajectory optimization across varying flight time and dynamic obstacle scenarios \cite{celestini2025generalizable}.
More recently, nonconvex trajectory generation based on natural language inputs with trajectory-level safety is demonstrated by the authors \cite{takubo2026semantic}. 
Despite these advances, existing learning-based frameworks still lack two critical capabilities for practical space mission design: (i) robust generalization across orbital regimes, and (ii) the ability to generate and explore Pareto-optimal trade-offs in multi-objective settings.
This paper addresses these limitations by developing a framework for generalizable nonconvex trajectory generation that approximates the true Pareto set.

\section{Background: Autonomous Rendezvous Transformer} \label{sec:background}

ART \cite{guffanti2024transformers} is a trajectory generation framework that leverages the sequence modeling and prediction capabilities of causal transformer models. 
Given a dataset of state and control trajectories, the core idea is to train the model via next-element prediction, learning to generate the optimal control inputs and corresponding states at each step.
A key enabler of ART for capable trajectory generation is the state representation, which is the tokenization of continuous trajectory information into sequences of states, controls, and associated performance metrics.
In ART, the tokenized trajectory representation is typically defined as
\small
\begin{align}
    \boldsymbol{\tau}_{1:N} = \{ \boldsymbol{g}_1, \boldsymbol{x}_1, r_1, c_1, \boldsymbol{o}_1,  \boldsymbol{u}_1, ..., \boldsymbol{g}_N, \boldsymbol{x}_N, r_N, c_N, \boldsymbol{o}_N,  \boldsymbol{u}_N \},
\end{align}
\normalsize
where $\boldsymbol{g}_k$ denotes the goal state, and $\boldsymbol{x}_k$ and $\boldsymbol{u}_k$ represent the state and control input at time step $k$, respectively. 
$N$ indicates the length of the time horizon, which has a constant time displacement $\Delta t$ for each time step. 
Additionally, $\boldsymbol{o}_k$ represents a scenario-specific observation vector at time step $k$ (e.g., information about obstacles in the scene).
This work considers two specific performance metrics, $ \{r_k, c_k\}$, where $r_k \in \mathbb{R}$ is the Reward-To-Go (RTG) and $c_k \in \mathbb{Z}_{\geq0}$ is the Constraint-To-Go (CTG). 
These metrics are defined as 
\begin{subequations}
\begin{alignat}{2}        
    r_{k} & = -\sum_{j=k}^{N} j_{j} (\boldsymbol{x}_j,\boldsymbol{u}_j,t_j), \quad 
    c_{k} & = \sum_{j=k}^{N} \boldsymbol{1}_{\mathcal{S}_{j}}(\boldsymbol{x}_j, \boldsymbol{u}_j), \ 
\end{alignat}
\end{subequations}
where $j_j$ is the instantaneous reward and $\boldsymbol{1}_{\mathcal{S}_{j}}(\cdot)$ is the indicator function for the feasible domain at the time step $j$. 
Essentially, both RTG and CTG measure the cumulative future optimality and feasibility of the trajectory starting from the current state.
A subsequence is defined as a tuple $[\boldsymbol{g}_k, \boldsymbol{x}_k, {r}_k, {c}_k, \boldsymbol{o}_k, \boldsymbol{u}_k]$, where each element type is embedded into a common input space of dimension $D=6$, corresponding to the number of distinct element types.
These embeddings are then projected into the transformer’s latent space.
This general framework allows tokenization of trajectory optimization problems with an arbitrary number of objectives and constraints. 

The workflow of ART comprises three phases: (i) dataset generation, (ii) training, and (iii) test-time inference.

\textbf{Dataset Generation:} The process begins by generating a set of scenario configurations through domain randomization of key problem parameters (e.g., boundary conditions, flight times). 
For each configuration, a corresponding trajectory optimization problem is solved to obtain reference solutions for training.
Diversity across parameters, operating conditions, and performance metrics is essential, as it ensures exposure to a broad solution landscape and supports generalization of the underlying dynamics and control strategies.

\textbf{Training:} After collecting a sufficient number of trajectories, the causal transformer model is trained to minimize a supervised loss. 
The final layer of the transformer branches into two task-specific heads, each implemented as a linear layer: one for state prediction and one for control input prediction. 
The loss function over a batch of $B$ trajectories with horizon $N$ is defined as
\small
\begin{align}
    \mathcal{L}(\boldsymbol{\tau}_{1:N}) = \sum_{b=1}^{B} \sum_{k=1}^{N} \left( \left\| \boldsymbol{x}_k^{(b)} - \hat{\boldsymbol{x}}_k^{(b)} \right\|_2^2 + \left\| \boldsymbol{u}_k^{(b)} - \hat{\boldsymbol{u}}_k^{(b)} \right\|_2^2 \right),
\end{align}
\normalsize
where $\hat{\boldsymbol{x}}_k$ and $\hat{\boldsymbol{u}}_k$ denote the model's predictions, and $\boldsymbol{x}_k$ and $\boldsymbol{u}_k$ represent the reference values from the dataset.

\begin{figure}
    \centering
    \includegraphics[width=0.9\linewidth]{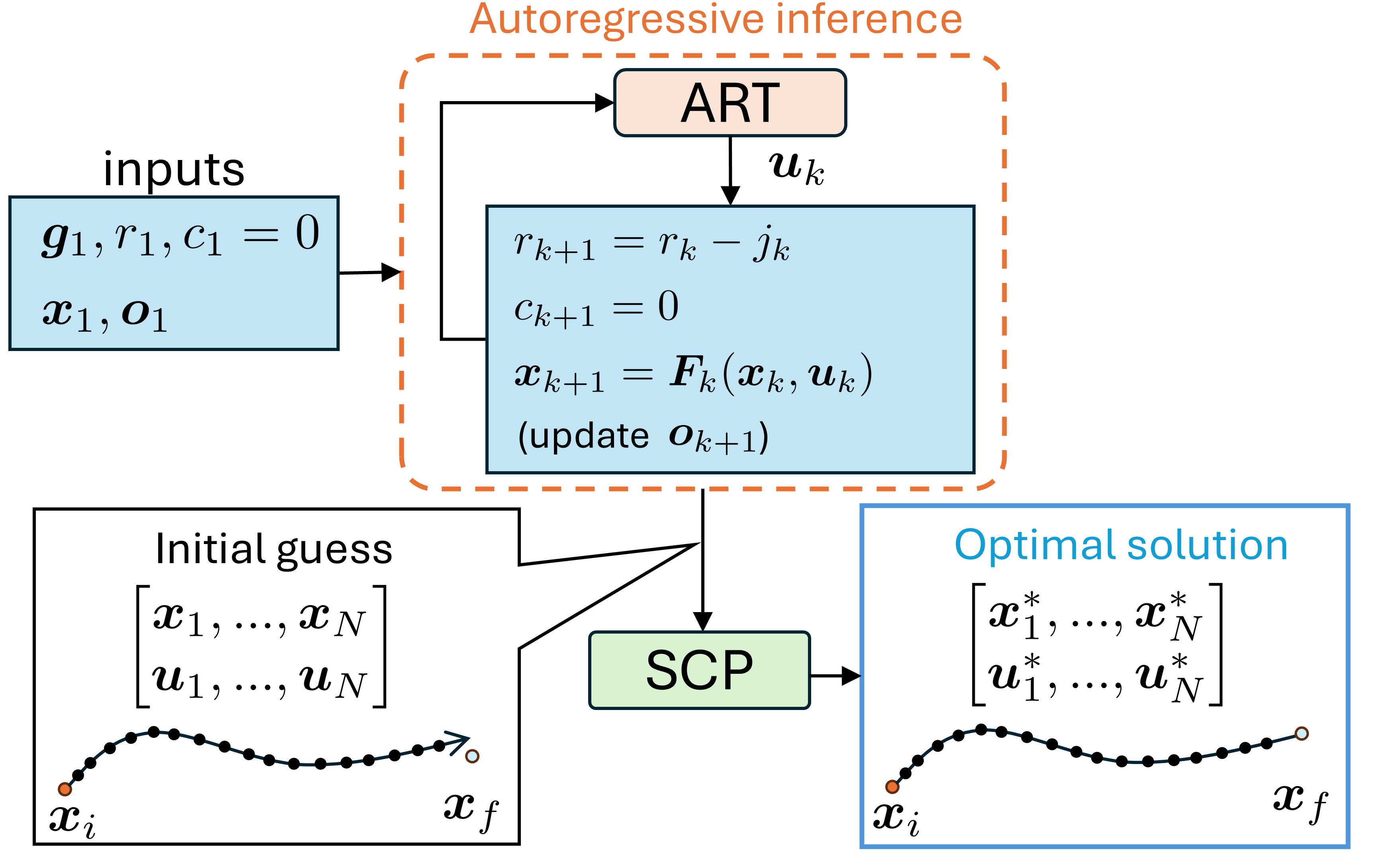}
    \caption{ART inference process. During autoregressive generation, ART only predicts the optimal control $\boldsymbol{u}_k$ while states are propagated via an available dynamics model.}
    \label{fig:art_inference}
\end{figure}

\textbf{Test-time inference:} After training, ART operates as an autoregressive trajectory generator, as illustrated in Figure~\ref{fig:art_inference}.
Given the initial input tokens (consisting of the desired terminal state, current state, and performance metrics, and observations), the transformer predicts the next control input. 
During inference, only the control input head is employed to generate the complete trajectory.
The predicted control input is then applied to the current state through an available dynamics model $\boldsymbol{F}(\boldsymbol{x}, \boldsymbol{u})$, which is also employed by the subsequent trajectory optimization algorithm. 
This model-in-the-loop rollout ensures that the predicted trajectories remain dynamically feasible.
The performance metrics are then updated at each step based on the control input and the propagated state.
Specifically, the RTG is initialized as a quantifiable lower bound on the optimal cost and updated as $r_{k+1} = r_k - j_k(\boldsymbol{x}_k,\boldsymbol{u}_k,t_k)$.
To incentivize constraint satisfaction, the CTG tokens are set to $c_1 = 0$ and kept at zero throughout inference (i.e., we condition on a zero-violation budget), which biases the policy toward feasible control input.
The cycle of control input prediction, dynamics propagation, and metric updating repeats until the terminal condition is met or the horizon is reached.
Note that this paper assumes the stationary goal state; therefore $\boldsymbol{g}_k$ is always set to the terminal state. 

\section{Methodology} \label{sec:methodology}

This section introduces the general class of multi-objective trajectory optimization problems addressed by the proposed approach, highlighting their formulation and integration with the ART framework. 

\subsection{Multi-Objective Optimal Control Problem}

Space mission design involves identifying feasible solutions while optimizing multiple, often competing, mission objectives. 
In the context of RPOD missions, a fundamental trade-off in preliminary trajectory design is balancing fuel consumption against flight time. 
Given this consideration, it is crucial for mission designers to characterize the Pareto front of the associated MO-OCP:
{\small
\begin{subequations}    
\begin{alignat}{2}
    \min_{\{\boldsymbol{x}_k\}, \{\boldsymbol{u}_k\}} \ & \{ \mathcal{J}_1, ..., \mathcal{J}_M \} \\
    \text{subject to} \quad & \boldsymbol{x}_{k+1} = \boldsymbol{F}_k (\boldsymbol{x}_{k}, \boldsymbol{u}_k) \quad && \forall k \in [1,N-1], \\
    & (\boldsymbol{x}_k, \boldsymbol{u}_k) \in \mathcal{S}_k 
    \quad && \forall k \in [1,N], \\
    & \mathcal{J}_l = \sum_{k=1}^N  j_{k,l}(\boldsymbol{x}_k,\boldsymbol{u}_k,t_k) \ && \forall l \in [1,M],  
\end{alignat}
\end{subequations}
}
where $\boldsymbol{x}_k = \boldsymbol{x}(t_k) \in \mathbb{R}^{n_x}$ and $\boldsymbol{u}_k = \boldsymbol{u}(t_k) \in \mathbb{R}^{n_u}$ represent the system state and the corresponding control input at time step $k$; 
$\boldsymbol{F}_k: \mathbb{R}^{n_x} \times \mathbb{R}^{n_u} \rightarrow \mathbb{R}^{n_x}$ denotes the discrete-time dynamics; and the set $\mathcal{S}_k \subseteq \mathbb{R}^{n_x} \times \mathbb{R}^{n_u}$ defines the admissible state-control pairs at time step $k$, encoding operational constraints. 
More generally, $\mathcal{S}_k$ may be expressed as a finite union of $n_c$ sets, $\mathcal{S}_k = \mathcal{S}_{k,1} \ \cup \ ... \ \cup \ \mathcal{S}_{k,n_c}$. 
Finally, ${\mathcal{J}_1, \ldots, \mathcal{J}_M}$ denotes the set of $M$ objective functions to be optimized.

There are two primary challenges in solving the above MO-OCP. 
First, multi-objective optimization problems are typically solved with computationally intensive metaheuristics \cite{Arya2024Stochastic}\cite{takubo2022robust}. 
A common alternative is to reformulate the problem as a single-objective optimization using either a weighted sum of objectives or a constraint-based scalarization \cite{vasile2011multi}, where all but one objective are converted into an inequality constraint with upper bounds. 
However, these reformulated problems can yield nondominated solution sets that are highly sensitive to minor variations in constraints, thereby limiting their robustness and generalizability.
Moreover, the problem is often nonconvex due to the presence of nonlinear dynamics and nonconvex feasible sets for the state and control variables.
SCP offers a promising solution strategy by iteratively solving convex approximations of the nonconvex problem.
This method leverages efficient convex solvers but necessitates a sufficiently accurate initial guess to ensure rapid convergence and to avoid poor local minima, as highlighted by Malyuta et al. \cite{malyuta_scp_2022}.

\subsection{ART for Efficient Pareto-Optimal Trajectory Generation}

Since ART was originally developed for single-objective, single-constraint problems, several extensions are introduced to enable trajectory generation over a multi-objective solution space. 
Specifically, the dimensions of the RTG and CTG are generalized to accommodate arbitrary numbers of rewards and constraints, such that $\boldsymbol{r}_k \in \mathbb{R}^{n_r}$ and $\boldsymbol{c}_k \in \mathbb{Z}_{\geq 0}^{n_c}$ denote the vectorized RTG and CTG, respectively. 
To achieve generalization across flight times and environmental parameters, the approach proposed by Celestini et al.~\cite{celestini2025generalizable} is adopted. 
Dataset generation employs domain randomization over three elements: (i) boundary conditions, (ii) flight times, and (iii) observation parameters defining the system dynamics. 
The randomization domain directly determines the domain over which the model generalizes, as induced by the seen data during the training.
In particular, flight time is controlled by a fixed step size $\Delta t$ together with a variable horizon length $N$, ensuring that the transformer’s learned dynamics evolve under a consistent transition interval.
This formulation simplifies training and improves efficiency at inference by allowing variable-length trajectories to be stacked and autoregressively generated within a single forward pass (cf. Figure~\ref{fig:art_md_summary}). 
Specifically, variable-horizon trajectories $\boldsymbol{\tau}_{1:N}^{(1:B)}$ are organized into a tensor of shape $[B, N, D]$, where $N$ denotes the maximum horizon length in the batch.
Shorter trajectories are padded to match the same length. 
At each step, the model autoregressively produces state–control input tokens, yielding $[B, D]$ outputs that are assembled into complete trajectories.

\section{Spacecraft Rendezvous in perturbed Keplerian orbits} \label{sec:case_study}
 
This paper focuses on mission design for on-orbit servicing operations targeting multiple resident space objects in Earth orbit. 
In this scenario, a servicing spacecraft (servicer) executes an approach maneuver to a predefined waypoint in the vicinity of a target spacecraft (target). 
During the preliminary trajectory design phase of such missions, it is essential to conduct rapid trade studies between flight time and total control cost across a diverse set of candidate targets, each characterized by distinct orbital elements and, consequently, unique relative dynamics.  
In addition, trajectory design must account for the presence of uncertainty in the system. 
Accordingly, this paper considers optimization of a nominal deterministic trajectory that explicitly models uncertainty effects and satisfies mission constraints in a probabilistic sense.

\subsection{Optimal Control Problem Formulation}
The trajectory optimization problem described above can be formulated as a discrete-time MO-OCP for chance-constrained, fault-tolerant rendezvous trajectory design under uncertainty~\cite{takubo2025towards} as
{\footnotesize
\begin{subequations}   
\begin{alignat}{2} \label{eq:opt_nonconvex}
    \text{given} \quad & \ \Delta t \\
    \min_{\{\boldsymbol{x}_k\}, \{\boldsymbol{u}_k\}} \ &  \left\{ \sum_{k=1}^{N} \|\boldsymbol{u}_k\|_2, \ N \right\} \\
    \text{subject to} \quad & \boldsymbol{x}_{k+1} = \Phi_k (\boldsymbol{x}_{k} + \Gamma_k \boldsymbol{u}_k), \quad  \forall k \in [1,N-1],  \label{eq:con_dyn} \\
    & \Pr( \boldsymbol{x}_{k0}^\top S_{ki} \boldsymbol{x}_{k0}  \geq 1) \geq 1-\Delta_k, \label{eq:con_ps_chance} \\
    & \qquad\qquad\qquad \qquad   \ \forall k \in [1,N-1], \forall i \in [1,N_s] \nonumber \\
    & \boldsymbol{x}_1 = \boldsymbol{x}_i, \quad  \boldsymbol{x}_{N} + \Gamma_N \boldsymbol{u}_N = \boldsymbol{x}_f, \label{eq:con_boundary} \\
    & \|\boldsymbol{u}_k\|_2 \leq u_{\text{max}} \label{eq:con_max_u}, \\ 
    \text{where} \quad & \Phi_{k} = \Phi(t_{k+1}, t_k), \quad \Phi_{ki} = \Phi(t_k + \tau_i, t_k), \\
    & \Psi_{ki} = \Psi(t_k+\tau_i), \\
    & S_{ki} = (\Psi_{ki} \Phi_{ki})^\top P_k \Psi_{ki} \Phi_{ki}, \\
    & \boldsymbol{x}_{k0} = \boldsymbol{x}_k+\Gamma_k\boldsymbol{u}_k.
\end{alignat}
\end{subequations}
}
The state $\boldsymbol{x} \in \mathbb{R}^6$ is expressed in Relative Orbital Elements (ROE) \cite{damico_phd_2010}\cite{delurgio2024closed}. 
Let the Keplerian orbital elements be $\boldsymbol{\alpha}_c = [a, e, i, \Omega, \omega, M]$, where $a$ is the semi-major axis, $e$ the eccentricity, $i$ the inclination, $\Omega$ the right ascension of the ascending node, $\omega$ the argument of periapsis, and $M$ the mean anomaly.
The ROE are a particular set of differences between the orbital elements of the target and the servicer (also called the deputy), denoted as $\boldsymbol{\alpha}$ and $\boldsymbol{\alpha}_d$, respectively. 
The objective is to simultaneously minimize the total control cost, modeled as a sum of impulsive maneuver magnitudes, and the terminal time represented by the discrete horizon length $N$, leading to a flgiht time of $N\Delta t$.  
Eq. \eqref{eq:con_dyn} describes the discrete-time ROE dynamics with impulsive velocity change $\boldsymbol{u} \in \mathbb{R}^3$ applied in the Radial/Tangential/Normal (RTN) frame. 
The State Transition Matrix (STM) $\Phi(t_{k+1}, t_k) \in \mathbb{R}^{6\times6}$ and control input matrix $\Gamma_k = \Gamma(t_k) \in \mathbb{R}^{6\times3}$ are expressed as an analytical function of the target's orbital elements but not the relative states themselves \cite{damico_phd_2010}\cite{delurgio2024closed}, which offers a convenient expression as a linear time-varying system.  
Additionally, Eqs.~\eqref{eq:con_boundary} and \eqref{eq:con_max_u} represent boundary conditions and maximum instantaneous control magnitude constraints, respectively.  
Finally, Eq. \eqref{eq:con_ps_chance} indicates the passive safety constraint, which ensures the deputy is outside of the ellipsoidal safety domains around the target not only for the entire duration of the controlled trajectory $t\in [0, t_f]$ but also for a predefined duration of the drift trajectories $\tau \in [0, \tau^s]$ after the sudden loss of control authority. 
The geometry of the safety domains at epoch $t_k$ is defined using the diagonal matrix $P_k \in \mathbb{R}^{6\times6}$, which encodes the semi-axes of the ellipsoid within the RTN frame. 
In this case study, two safety domains are considered: the Approach Ellipsoid (AE) and the Keep-Out Zone (KOZ), both modeled as spheres. 
At the switching epoch, $t_\text{switch}$, the active safety domain transitions from the larger AE to the smaller KOZ, reflecting the practical scenario in which the allowable safety region contracts as the chaser spacecraft closes in on the target.
Matrix $\Psi_k = \Psi (t_k) \in \mathbb{R}^{6\times6} $ denotes the first-order linear mapping between the ROE and the matrix expression of the relative position and velocity vector in the Cartesian coordinates within the RTN frame (denoted as $\boldsymbol{x}_{\text{RTN}}$) \cite{damico_phd_2010}.
The passive safety constraint is formulated as a chance constraint to account for the inherent uncertainty that comes from navigation, control execution error, and unmodeled acceleration. 
At each epoch, the risk factor (i.e., failure probability) $\Delta_k$ is chosen as a parameter, and the probabilistic satisfaction of the constraint is enforced. 
Assuming Gaussian state uncertainty, the drifted mean state and covariance propagated from the time immediately after the control epoch $t_k$ over the drift duration $\tau_i$, denoted by $\boldsymbol{x}_{ki} = \mathbb{E}[\boldsymbol{x}(t_k; \tau_i)]$ and $\Sigma_{ki} = \Sigma(t_k; \tau_i) = \mathbb{E}[\boldsymbol{x}(t_k; \tau_i)\boldsymbol{x}(t_k; \tau_i)^\top]$, can be written as 
{\footnotesize
\begin{subequations} \label{eq:con_ps_det}
\begin{align}
    & 1 - \boldsymbol{x}_{k0}^\top S_{ki} \boldsymbol{x}_{k0} + q(\Delta_k) \sqrt{\boldsymbol{g}_{ki}^\top \Sigma_{ki} \boldsymbol{g}_{ki}}  \leq 0 \\
    & \boldsymbol{g}_{ki} = -2 \Phi_{ki}^{-\top}(S_{ki}   \boldsymbol{x}_{k0}), \\
    & \Sigma_{k(i+1)} = \Phi_{ki} \Sigma_{ki} \Phi_{ki}^\top + Q_{ki}, \quad   \Phi_{ki} = \Phi(t_k+\tau_{i+1}, t_k+\tau_{i}), \\
    & \Sigma_{k0} = \Sigma_{\text{nav}}(\boldsymbol{x}_{k}) + \Gamma_k \Sigma_{\text{exe}} (\boldsymbol{u}_k) \Gamma_k^\top, \label{eq:con_cov_ini}
\end{align}
\end{subequations}
}
where $q(\cdot)$ is the inverse of the normal cumulative distribution function. 
Three uncertainty sources are considered in this case study: unmodeled acceleration, navigation uncertainty, and control execution error. 
$Q_{ki}$ is the discrete-time process noise covariance, representing the white-noise unmodeled acceleration.  
Eq. \eqref{eq:con_cov_ini} indicates that the covariance immediately after the loss of control has contributions from both navigation uncertainty and the control execution error. 
A simplified navigation model is introduced to model the navigation covariance as a function of the range between the target and the servicer as \cite{takubo2025towards}
\begin{equation} \label{eq:artms_nav_general}
    \Sigma_{\text{nav}}(\boldsymbol{x}_k) = \rho(\boldsymbol{x}_k) \boldsymbol{s} \boldsymbol{s}^\top,
\end{equation}
where $\rho(\boldsymbol{x}_k)$ is the range between the target and the servicer, and $\boldsymbol{s} \in {\mathbb{R}}_{\geq0}^{6}$ is a user-defined parameter.
Gates model \cite{gates1963simplified} is employed for the execution error covariance, denoted as $\Sigma_{\text{exe}}$, which is a function of the nominal control input. 
In particular, this is expressed in the frame based on the thrust direction, denoted here by the superscript $V$ as 
\begin{equation} \label{eq:gates_model}
\resizebox{0.9\linewidth}{!}{$
    \Sigma_{\text{exe}}^V(\boldsymbol{u}) = {\mathrm{Diag}}\{\sigma_{m1}^2 + \sigma_{m2}^2 \| \boldsymbol{u} \|_2 ^2 , \sigma_{p1}^2  + \sigma_{p2}^2 \| \boldsymbol{u} \|_2^2, \sigma_{p1}^2 + \sigma_{p2}^2 \| \boldsymbol{u} \|_2^2\} ,
$}
\end{equation} 
where $\{\sigma_{m1}, \sigma_{m2},\sigma_{p1},\sigma_{p2}\} \in {\mathbb{R}}_{\geq0}^{4}$ are model parameters. 
Evidently, this constraint is nonconvex with respect to $(\boldsymbol{x}_k, \boldsymbol{u}_k)$, which is convexified at each iteration of the SCP. 
Furthermore, the number of constraints is reduced from $NN_s$ to $N$ by only evaluating the worst epoch $i^*\in [1,N_s]$ within each drift trajectory in the previous iteration \cite{guffanti2022phd}.

\subsection{Relative Orbital Elements}
Two sets of ROEs are used for the state representation in this paper. 
First, quasi-nonsingular ROE (qnsROE) is one of the most well-known ROE sets, defined as \cite{damico_phd_2010}
{\footnotesize
\begin{equation}
\boldsymbol{x}_{\text{qns}} = \left[\begin{array}{c}
\delta a \\ \delta \lambda \\ \delta e_x \\ \delta e_y \\ \delta i_x \\ \delta i_y
\end{array}\right]
=
\left[\begin{array}{c}
\left(a_d-a\right) / a \\
(M_d+\omega_d) - (M+\omega) +\left(\Omega_d-\Omega\right) \cos i \\
e_d \cos \omega_d - e  \cos \omega  \\
e_d \sin \omega_d - e  \sin \omega  \\
i_d-i \\
\left(\Omega_d-\Omega\right) \sin i
\end{array}\right].
\end{equation}
}
The qnsROE is well-suited for describing relative motion in near-circular orbits, as it provides an elegant geometric description of bounded relative motion \cite{damico_phd_2010}.
In this paper, an STM that accounts for the secular $J_2$-effect \cite{koenig2017new} is adopted. 
The first-order control input matrix is provided in D'Amico \cite{damico_phd_2010}. 
Both the STM and control input matrices are linear time-varying, determined by the target’s orbital elements.
In elliptic orbits, the qns-ROE no longer preserves the simple geometric correspondence with the osculating relative motion.
To address this issue, this study adopts the elliptic ROE (eROE), which is first-order equivalent to the integration constants of Yamanaka–Ankersen STM \cite{yamanaka2002new}\cite{guffanti2022phd}. 
The eROE is defined as \cite{delurgio2024closed}
{\footnotesize
\begin{equation}
\boldsymbol{x}_{\text{e}}
= \left[\begin{array}{c}
\delta a \\ \delta \lambda \\ \delta e_x \\ \delta e_y \\ \delta i_x \\ \delta i_y
\end{array}\right]
=
\left[\begin{array}{c}
\eta^2\left(a_d-a\right) / a \\
\frac{1}{\eta}\left(M_d-M\right)+\eta^2\left[\omega_d-\omega+\left(\Omega_d-\Omega\right) \cos i \right] \\
\left(e_d-e\right) \cos \omega + \frac{e}{\eta}\left(M_d-M\right) \sin \omega \\
\left(e_d-e\right) \sin \omega - \frac{e}{\eta}\left(M_d-M\right) \cos \omega \\
\eta^2\left(i_d-i\right) \\
\eta^2\left(\Omega_d-\Omega\right) \sin i
\end{array}\right],
\end{equation}
}
\noindent where $\eta = \sqrt{1-e^2}$. 
Similar to the qnsROE, the mean $J_2$-perturbed STM and the control input matrix are available in Refs. \cite{delurgio2024closed}\cite{guffanti2022phd}.

\subsection{Sequential Convex Programming}
In this study, \texttt{SCVx$^\ast$} \cite{oguri2023successive} is employed as the SCP algorithm. 
This method ensures the feasibility of convexified subproblems by incorporating violations of the linearized nonconvex constraints into a penalty, formulating an augmented Lagrangian \cite{mao2016successive}. 
The trust region and penalty weights are automatically updated until convergence is achieved, with the hyperparameters provided in Table \ref{tab_appendix:scp_param} in Appendix. 
During the SCP routine, all variables are normalized with respect to the initial guess, thereby scaling each element to the range [-1,1] to enhance numerical stability \cite{malyuta2023fast}.
Specifically, the violations of the linearized passive safety constraint (cf. Eq.~\eqref{eq:con_ps_det}) and the terminal constraint (cf. Eq.~\eqref{eq:con_boundary}) are incorporated as penalty terms. 
This treatment is necessary despite the terminal constraint being convex, because ART's initial guess does not strictly satisfy the boundary conditions; by softening this constraint, the initial infeasibility of the convex subproblem is avoided.

Although theoretical convergence guarantees are discussed in various SCP algorithms \cite{oguri2023successive}\cite{mao2016successive}\cite{bonalli_2019_gusto}, non-converging scenarios are still observed in practice. 
Such cases typically occur due to (i) the absence of a feasible domain in the original nonconvex problem, (ii) numerical instabilities (particularly when the augmented Lagrangian penalty weight becomes excessively large), (iii) failure to converge within the prescribed maximum number of iterations. 
Furthermore, the aforementioned constraint reduction scheme for the nonconvex passive safety constraints also makes the linearization inaccurate and contributes to the non-converging behavior. 

\section{Results and Analysis}  \label{sec:results}

\begin{table*}[t!]
\caption{Domain Randomization used for the training dataset generation. $\textsf{linspace}(x_{\text{min}}, x_{\text{max}}, n)$ represents an $n$-element vector of evenly spaced points between $x_{\text{min}}$ and $x_{\text{max}}$, including the two boundary values.}
\label{tab:prob_setup}
\centering
\renewcommand{\arraystretch}{1}
{\footnotesize
\begin{tabular}{cccc}
\thickhline
 &  & LEO & Elliptic orbit \\
\hline
\multirow{6}{*}{$\alpha_c (t_1)$} 
 & $a $    & $6738.14 \ [\text{km}]$ & $\textsf{linspace}(6.5\cdot 10^3, 11\cdot 10^3, 10) \ [\text{km}]$ \\
 & $e$  & $5.58 \cdot 10^{-4}$ &  $\textsf{linspace}(10^{-3}, 0.35, 10)$ \\
 & $i$      & $51.641  \ [^\circ]$ &  $\textsf{linspace}(1 \cdot 10^{-3}, \pi-1 \cdot 10^{-3}, 10) \ [\text{rad}]$ \\
 & $\Omega $ & $301.037 \  [^\circ]$ & $\textsf{linspace}(0, 9\pi/5, 10) \ [\text{rad}]$ \\
 & $\omega $  & $26.18 \ [^\circ]$ & $\textsf{linspace}(0, 9\pi/5, 10) \ [\text{rad}]$ \\
 & $M $       & $68.23 \ [^\circ]$ & $\textsf{linspace}(0, 9\pi/5, 10) \ [\text{rad}]$ \\
\hline
\multicolumn{2}{c}{Time} 
 & \makecell{$\Delta t = 171.4\, \text{s}$,  $N = [30, 35, 40, 45, 50]$ \\ (when $N=50$, $t_f = 140$ min.)}   & \makecell{$\Delta t = 440.8\, \text{s}$,  $N = [35, 38, 41, 44, 47, 50]$ \\ (when $N=50$, $t_f = 360$ min.)}  \\
\hline
\multicolumn{2}{c}{State representation} & qnsROE & eROE \\ 
\hline 
\multicolumn{2}{c}{Initial state } & $x_{\text{RTN}} = \begin{bmatrix}
   -4.0 \pm \textsf{linspace}(-1.4, 1.4, 1000)\\
   -17.5 \pm \textsf{linspace}(-3.0, 3.0, 1000) \\
   0 \\
   0 \\
   6.849 \pm  \textsf{linspace}(-1.0, 1.0, 100) \\
   0 
\end{bmatrix}$ [km, m/s] & \makecell{ $\Delta a \sim \textsf{linspace}(-20.0, -7.0, 10) $ [km] \\ 
 $\Delta e \sim \textsf{linspace}(-3 \cdot 10^{-4}, 3 \cdot 10^{-4}, 10) $ \\ 
 $\Delta i \sim \textsf{linspace}(-6 \cdot 10^{-4}, 6 \cdot 10^{-4}, 10) \ [\text{rad}] $ \\ 
 $\Delta \Omega \sim \textsf{linspace}(-3 \cdot 10^{-4}, 3 \cdot 10^{-4}, 10) \ [\text{rad}]$ \\ 
 $\Delta \omega \sim \textsf{linspace}(-3 \cdot 10^{-3}, 3 \cdot 10^{-3}, 10) \ [\text{rad}]$ \\ 
 $\Delta M = -0.03 \ [\text{rad}] $ } \\
\hline
\multicolumn{2}{c}{Terminal state [km, m/s]} & 
$\boldsymbol{x}_{\text{RTN}} = \begin{cases}
    [ \pm 0.75, 0,0, 0,0,0 ]\\
    [ 0, \pm 0.75, 0,0,0,0 ]
\end{cases}$ & 
$\boldsymbol{x}_{\text{RTN}} = [ 0,1.5,0, 0,0,0 ] $  \\
\thickhline 
\multicolumn{2}{c}{$u_{\max} \ [\text{m/s}]$} & 5.0 & 7.0 \\
\hline
\multicolumn{2}{c}{$(\tau^s \ [\text{min.}], N_s)$} & (91.7 ($\sim$ 1 orbit), 30)  & (360, 30) \\
\multicolumn{2}{c}{$(r_{\text{AE}} \ [\text{km}], r_{\text{KOZ}} \ [\text{km}], t_{\text{switch}} [\text{min.}])$} & ( $1.0, 0.5, 77 $ ) & ( $6.0, 1.0, 198 $ )\\
\thickhline
\end{tabular}
}
\end{table*}

\begin{figure*}[ht!]
    \centering
     \begin{subfigure}[ht!]{0.4\textwidth}
         \centering
         \includegraphics[width=\textwidth]{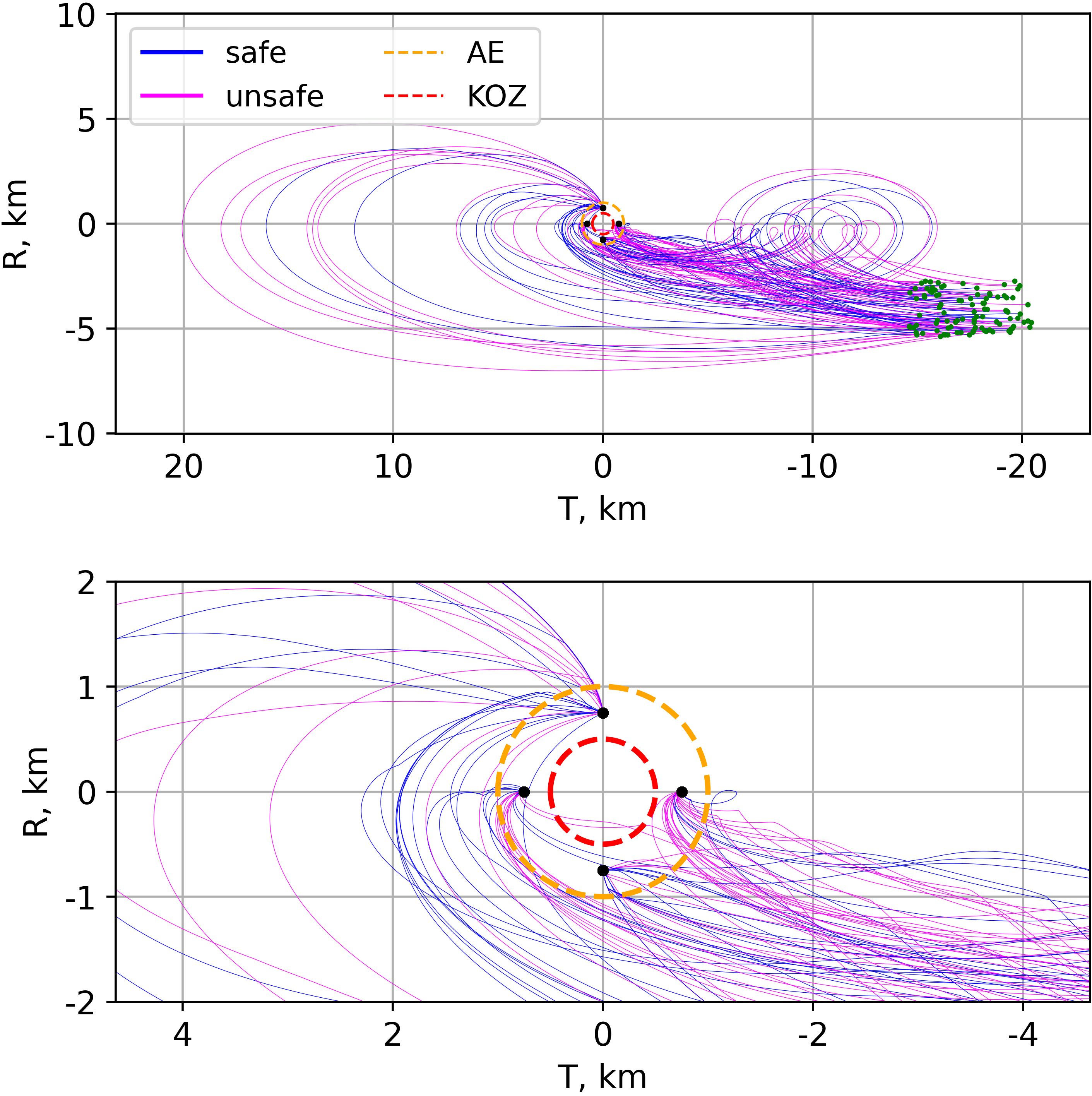}
         \caption{LEO scenario.}
         \label{fig:leo_traj_RT}
     \end{subfigure}
     \begin{subfigure}[ht!]{0.4\textwidth}
         \centering
         \includegraphics[width=\textwidth]{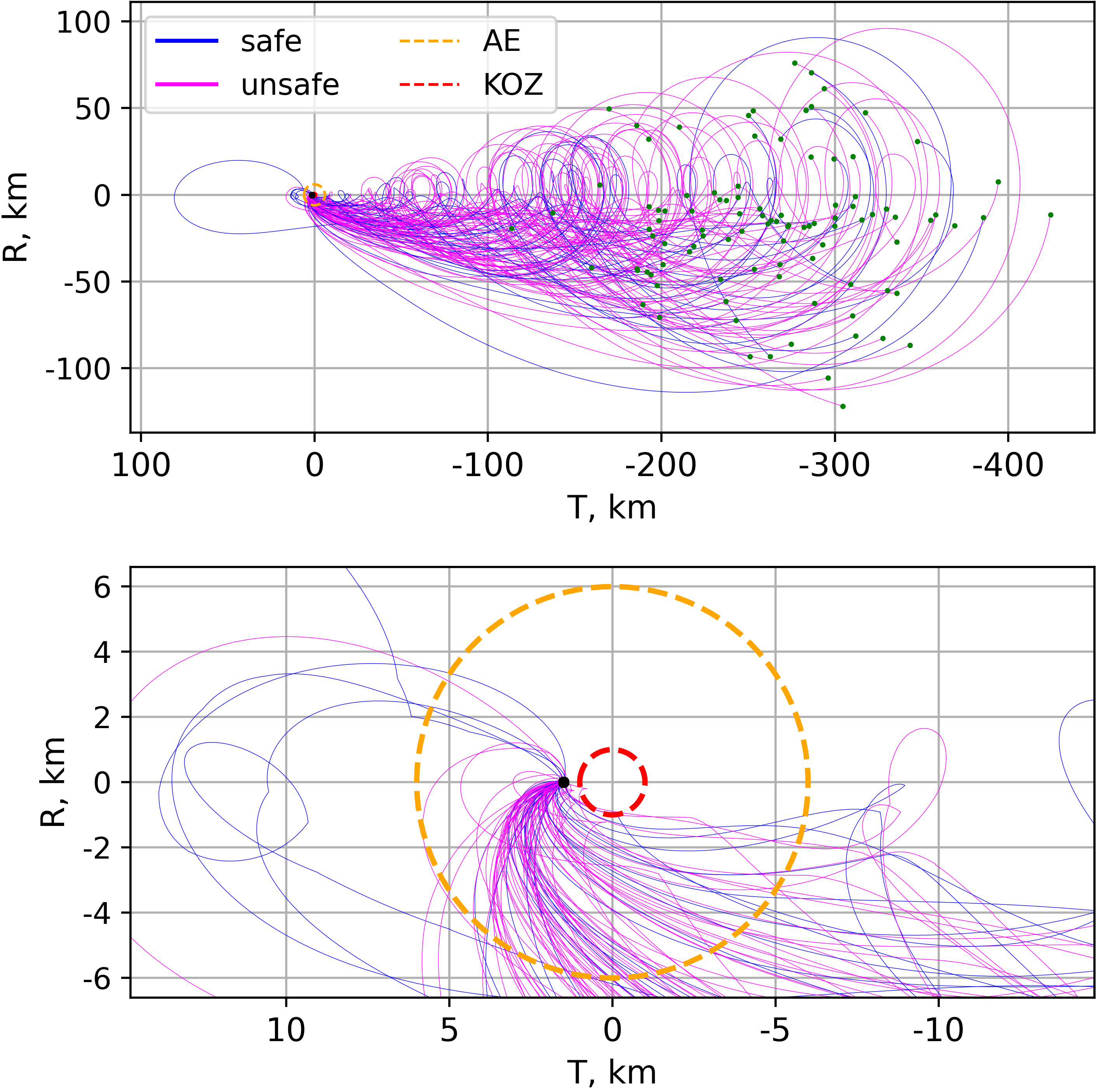}
         \caption{Generalized elliptic orbit scenario.}
         \label{fig:kep_traj_RT}
    \end{subfigure}
    \caption{Sampled 100 trajectories from the dataset for two rendezvous scenarios. The lower panels provide zoomed-in views of the corresponding upper panels. Orange and red dotted circles indicate the two-dimensional projection of the AE and KOZ, respectively.}
    \label{fig:traj_data}
\end{figure*}

This section evaluates the performance of the ART framework through a series of numerical experiments\footnote{The numerical experiments are implemented in Python and executed on a workstation equipped with an AMD Ryzen 9 7950X CPU and an NVIDIA GeForce RTX 4090 GPU.}. 
In the context of AI-driven space mission design, the proposed framework aims to serve as a flexible, accurate, and efficient generator of trajectories and mission profiles for a broad spectrum of RPOD missions. 
Accordingly, the study is guided by two central Research Questions (RQs): \textbf{(RQ1) Can ART reliably generate high-quality trajectories across varying flight times and diverse orbital dynamics?} and \textbf{(RQ2) Can ART serve as an accurate surrogate for the Pareto front?} 
The following discussion addresses these questions by demonstrating the capabilities of the proposed transformer-based framework for mission design.

\subsection{Scenarios and Benchmarks}

The training dataset is constructed by solving a batch of single-objective SCP problems, each initialized from a convex solution with a selected time horizon $N$ that omits the passive safety constraint in Eq. \eqref{eq:con_ps_chance} and minimizes the total control cost.
Table \ref{tab:prob_setup} summarizes the scenario setup and the domain randomization used to generate the training dataset.
Specifically, during the dataset generation, the flight time is randomly drawn from a set of evenly spaced values between two extremes.
In total, 80,000 problem instances are generated.
For each instance, both the convex relaxation and the corresponding SCP solution initialized with the convex solution are computed and stored. 
Only those trajectories for which the SCP successfully converges are retained, ensuring that the dataset consists of an even distribution of convex–nonconvex solutions, where the nonconvex solution is feasible with respect to the original optimization problem.
Additionally, during inference, the model initializes the RTG value using the objective of the convex relaxation, which is the lower bound under the time discretization based on the fixed time step $\Delta t$.

To serve as a flexible mission design tool, trajectories must be generated reliably and efficiently across a range of scenario configurations.  
With this in mind, two representative rendezvous mission scenarios are considered.  

\subsubsection{Mid-range low Earth orbit rendezvous}  
The first scenario builds on prior work~\cite{takubo2025towards} concerning the design of mid-range approach maneuvers in a Low Earth Orbit (LEO), where the target orbit is fixed while the flight time is varied. 
Because the target orbit is fixed, no observation vector is provided when training the model. 
The boundary conditions are imposed on the Radial–Tangential (RT) plane (cf. Table \ref{tab:prob_setup}).
As a result, the obtained rendezvous trajectories lie predominantly within the RT-plane, with the servicer exploiting drift in $\delta \lambda$ through semimajor axis differences to approach the target.
This profile is representative of crewed RPOD missions~\cite{dsouza2007orion}.  
However, such configurations lack inherent passive safety in relative motion~\cite{damico_phd_2010}, motivating the need for nonconvex trajectory optimization.  
The problem setup is summarized in Table~\ref{tab:prob_setup}, and Figure~\ref{fig:leo_traj_RT} illustrates 100 sample trajectories from the training dataset projected onto the RT-plane.  
Parameters related to execution error, navigation modeling, unmodeled acceleration, and risk factors are summarized in Table~\ref{tab_appendix:uncertainty} in Appendix.  

\subsubsection{Far-range elliptic orbit rendezvous}  
The second scenario addresses far-range rendezvous in various elliptic orbits, where both the flight time and the orbital elements of the target are modulated.  
To provide the dynamics context, the target's orbital elements are embedded as the observation vector at each time step as $\boldsymbol{o}_k =
[a, e, \sin i, \cos i, \sin \Omega, \cos \Omega, \sin \omega, \cos \omega, \sin M, \cos M ] \in \mathbb{R}^{10}$. 
This avoids discontinuities arising from angular variables jumping from $2\pi$ to $0$, and helps the transformer model better understand the context. 
Because the boundary conditions are not bound on the RT-plane (cf. Table \ref{tab:prob_setup}), the generated trajectories via SCP exhibit natural motion in the normal direction. 
However, owing to the larger eccentricity and the placement of the terminal waypoint near the boundary of the KOZ on the RT-plane, the resulting optimization problem remains concentrated in the vicinity of the RT-plane and is therefore particularly challenging.
Note that the terminal waypoint defined in the RTN frame may not satisfy the passive safety constraint, depending on the target's orbital elements.
Therefore, in this case study, the passive safety constraint is intentionally relaxed only at the final step to avoid the problem becoming infeasible due to the conflict between the terminal state equality and the passive safety constraint.
The process of generating terminal (or mid-course) waypoints based on the target's orbit is beyond the scope of this paper and is left for future work.
Moreover, as eccentricity increases, many orbits intersect the Earth's surface in practice.
This is, however, an intentional choice made to enhance generalizability across different relative dynamics scenarios.
As with the LEO case, the problem setup is summarized in Table~\ref{tab:prob_setup}, and Figure~\ref{fig:kep_traj_RT} presents 100 representative trajectories projected onto the RT-plane.  

\begin{figure*}[ht!]
    \centering
     \begin{subfigure}[ht!]{\textwidth}
         \centering
         \includegraphics[width=\textwidth]{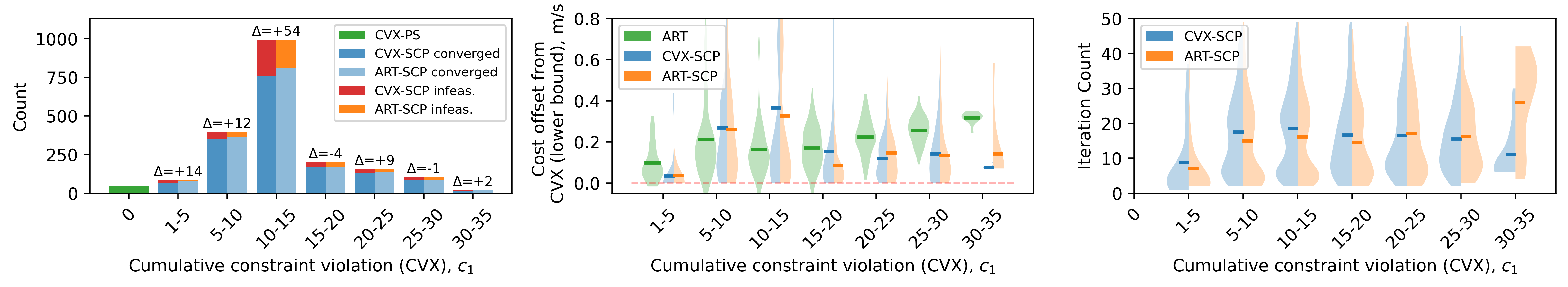}
         \caption{Generalization across flight time with the LEO scenario.}
         \label{fig:leo_analysis}
     \end{subfigure}
     \vspace{-0.5em} 
     \begin{subfigure}[ht!]{\textwidth}
         \centering
         \includegraphics[width=\textwidth]{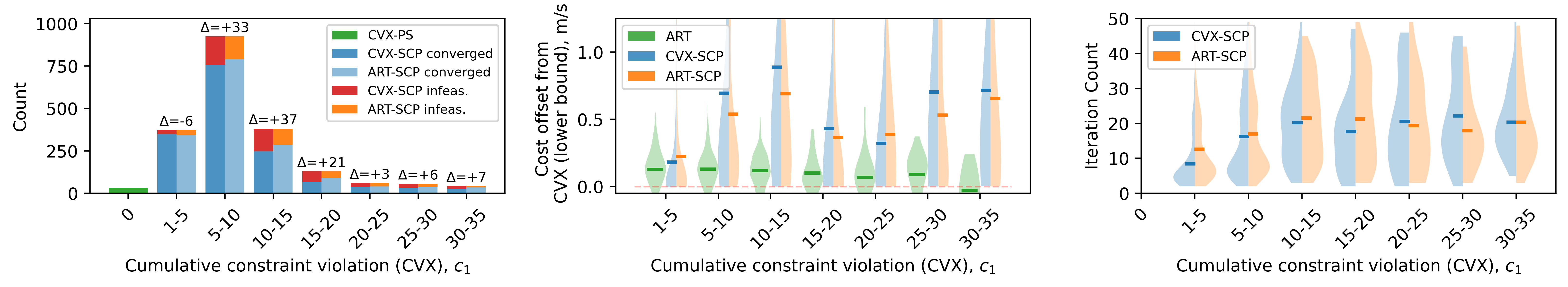}
         \caption{Generalization across dynamics and flight time with the elliptic orbit scenario. The test scenarios are not seen by the model during the training, yet are sampled from the same distribution. }
         \label{fig:kep_analysis_test}
    \end{subfigure}
    \vspace{-0.5em} 
    \begin{subfigure}[ht!]{\textwidth}
         \centering
         \includegraphics[width=\textwidth]{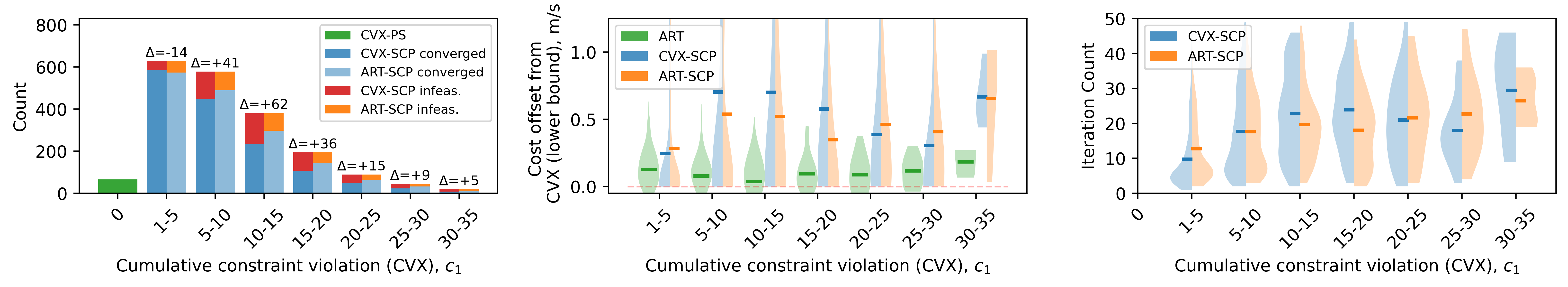}
         \caption{Generalization across the dynamics and flight time with the elliptic orbit scenario. The test scenarios are taken from SpaceTrack's medium-sized debris catalog (cf. Figure~\ref{fig:debris_oe}).  }
         \label{fig:kep_analysis_debris}
    \end{subfigure}
    \caption{Performance analysis of the SCP with ART initialization compared to that with CVX initialization, tested for 2,000 test scenarios. The statistics are binned across the cumulative constraint violation of the CVX solution.}
    \label{fig:ws-analysis}
\end{figure*}

The performance of ART is assessed by comparing the resulting SCP solutions initialized with ART-generated trajectories against those obtained using a conventional convex initialization.
In particular, a distinct set of problem instances, separate from those used in the training dataset and referred to as test scenarios, is employed for the analysis. 
The design of the test scenarios for each analysis is described in the following sections.
For simplicity, the following notation is used in the remainder of the paper:
\begin{itemize}
    \item CVX: solution of the convex minimum-fuel problem without passive safety
    \item CVX-SCP: passively safe solution obtained via SCP, initialized with CVX trajectory
    \item ART: trajectory generated by the ART inference 
    \item ART-SCP: passively safe solution obtained via SCP, initialized with ART trajectory. 
\end{itemize}

\subsection{RQ1: Generalization across Flight Time and Dynamics}

The generalizability of the ART inference is evaluated based on the quality of the generated trajectories as warm starts for the subsequent SCP.  
To this end, 2,000 test scenarios are generated, and the performance of CVX-SCP and ART-SCP is compared across three statistical metrics: (i) the success rate of convergence in SCP, (ii) the optimality of the resulting objective values, and (iii) the number of iterations required for SCP convergence.  
It is important to note that the baseline CVX solutions display a wide variation in cumulative trajectory-level constraint violations (i.e., $c_1$), with some trajectories even satisfying the nonconvex constraints incidentally.
These cumulative constraint violations provide an indicator of the difficulty in converging to feasible solutions.  
Accordingly, the problems are first grouped into bins based on the cumulative constraint violation of the CVX solutions, after which the statistical analyses are performed.  
Three experiments are conducted to evaluate ART's generalizability, as summarized in Figure \ref{fig:ws-analysis}.
For each experiment, results are presented in three forms: (i) histograms indicating whether the SCPs converged or not, (ii) violin plots of optimality, measured as the deviation from the CVX solution (lower bound), and (iii) violin plots of the iteration counts required for the SCP to converge. These visualizations are organized according to the binned groups of cumulative constraint violation in the CVX solutions.

First, generalization over flight time \cite{celestini2025generalizable} is evaluated in the context of the LEO rendezvous scenario \cite{takubo2025towards}.
Figure~\ref{fig:leo_analysis} presents the statistical analysis of the ART-generated trajectories. 
For this analysis, the test scenarios consist of previously unseen instances sampled from the same distribution as the training dataset.
The results show that initialization with the ART-generated trajectories leads to a higher success rate of the convergence in the SCP across a wide range of constraint violations in the CVX solutions. 
Meanwhile, the optimality of the converged solution and the iteration counts to convergence are statistically improved or equivalent to the CVX-warmstarting.  
Overall, the results demonstrate that ART learns an effective control policy across varying flight times with reliability comparable to that of traditional convex-based methods. 
Moreover, it provides effective initial guesses that improve success rates, enhance optimality, and accelerate convergence in the subsequent SCP.

Next, generalization across both orbital domains and flight times is evaluated by varying the targets’ orbital elements, thereby inducing different underlying relative dynamics. The performance of the trained model under these conditions is summarized in Figure~\ref{fig:kep_analysis_test}.
The test scenarios are drawn from the same distribution as the training dataset. 
The results present a similar trend to that of the LEO rendezvous scenario, where the success rate of the SCP's convergence using ART initialization outperforms the CVX initialization across most of the $c_1$ bins, while both the achieved optimality and the number of iterations required for convergence remain equivalent or superior when using ART-generated initial guesses. 

\begin{figure}[t!]
    \vspace{0.5em}
    \centering
    \includegraphics[width=\columnwidth]{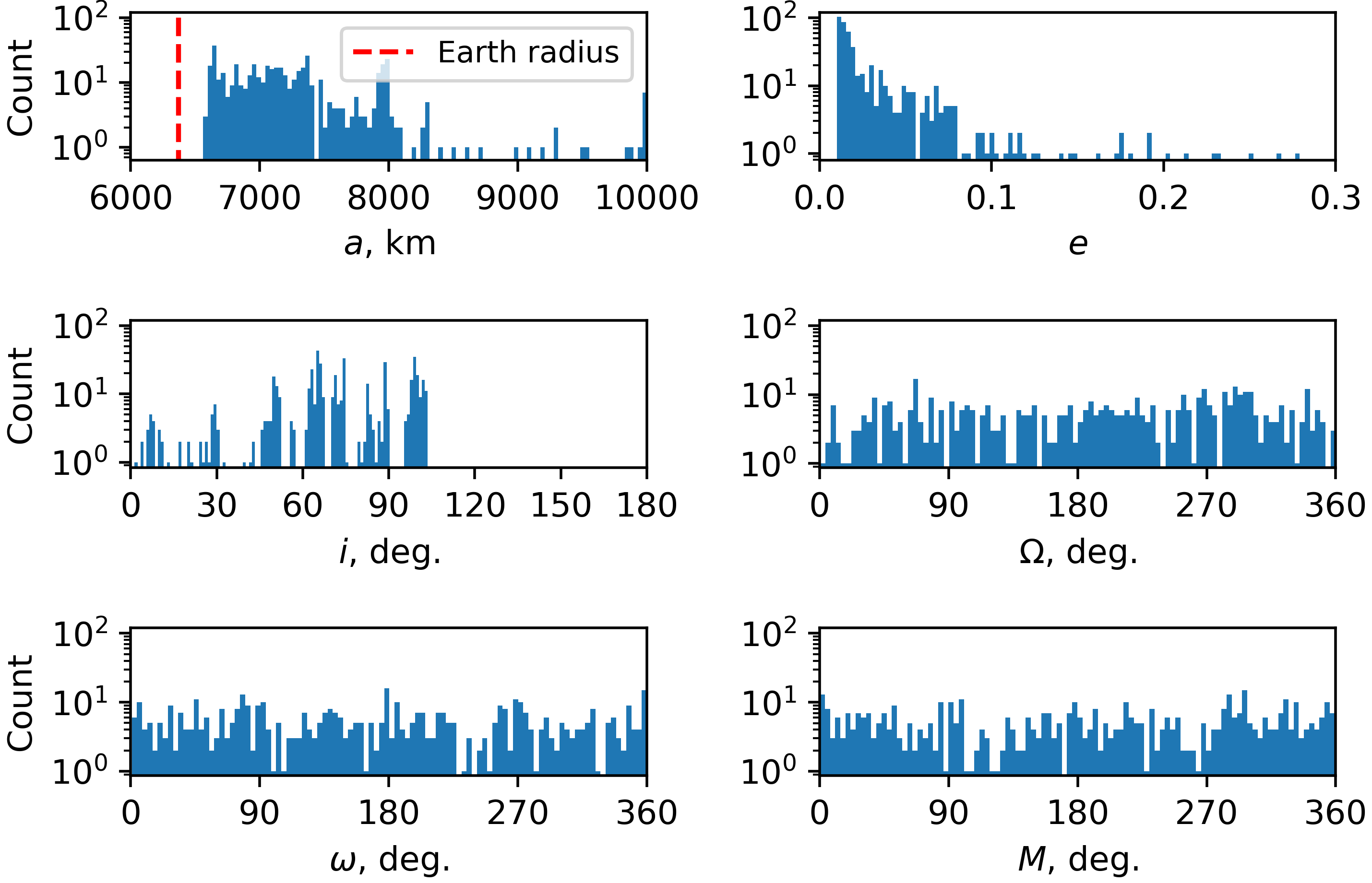}
    \caption{Distributions of the orbital elements of 500 sampled debris used for the test scenarios.}
    \label{fig:debris_oe}
\end{figure}

As a demonstration of generalization, the model trained with the elliptic orbit rendezvous dataset with varying flight times and underlying dynamics is further validated using the test scenarios based on the SpaceTrack orbital debris database\footnote{\url{https://www.space-track.org}. Data acquisition date: Jan. 4, 2026}. 
From a catalog of 5,000 medium-sized debris objects, 500 orbital elements are randomly selected within the ranges of semimajor axis ($6,500<a<10,000$ km) and eccentricity ($1\cdot10^{-3}<e<0.3$).
To promote diversity in the sampled orbital elements, targets are chosen in proportion to the logarithmic density across semimajor axis and eccentricity within each $(a,e)$ grid cell, thereby compensating for the strong concentration of debris in near-circular LEO. 
The resulting distribution of orbital elements of the target is shown in Figure~\ref{fig:debris_oe}.
Although the sampled orbital elements lie within the randomized training ranges, their values do not coincide with the discretized grid points, creating more challenging test scenarios for the model. 
The corresponding initial and boundary states are then generated by random sampling according to the specifications in Table~\ref{tab:prob_setup}.
Overall, a total of 2,000 test scenarios are generated and evaluated. 

The results of the warm-starting performance evaluation are summarized in Figure~\ref{fig:kep_analysis_debris}. 
Favorable results in terms of success rate, optimality, and iteration count are consistently observed across the cumulative constraint violation of the CVX solutions, particularly for cases with $c_1 > 5$. 
Although the training dataset may contain up to $10^6$ grid points of the chief orbital elements (cf. Table~\ref{tab:prob_setup}), the test scenarios considered in this experiment still require interpolation over a large six-dimensional continuous space. 
This observation supports the capability of the trained model to generalize to previously unseen parameter values within the range of the training dataset.

\subsection{RQ2: Prediction Capability of Mission Design Objectives}

\begin{figure}[t!]
     \centering
     \vspace{0.5em}
     \begin{subfigure}[th!]{\columnwidth}
         \centering
         \includegraphics[width=\columnwidth]{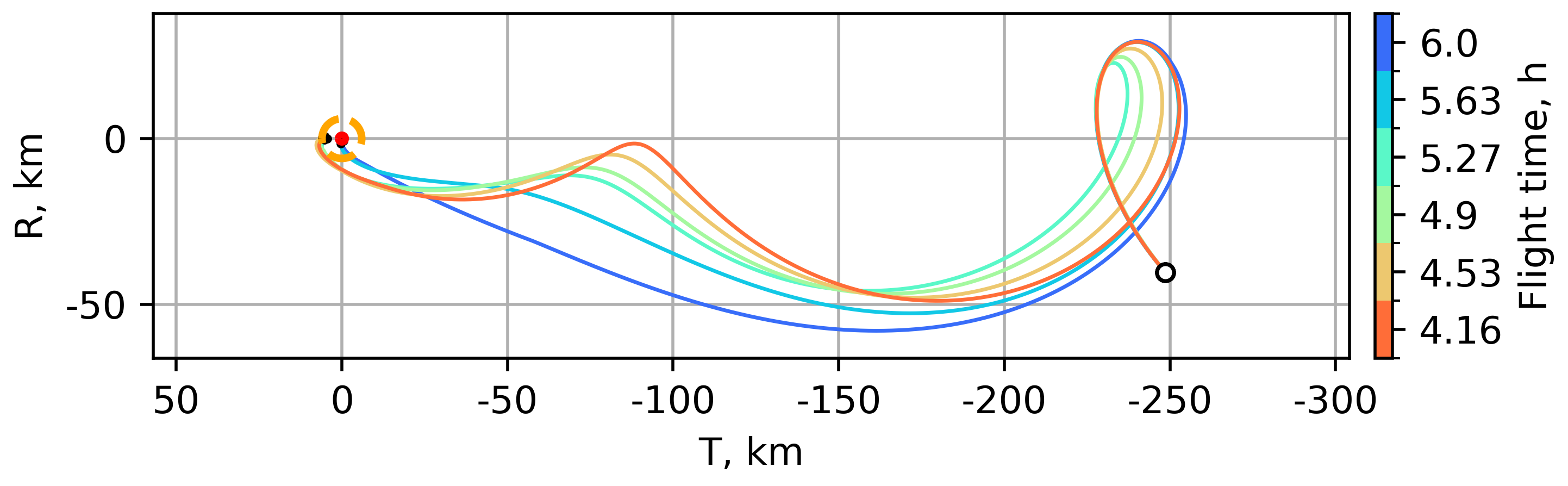}
         \caption{Batch generation of near-optimal and passively-safe trajectories from the same initial state and target waypoint, with varying flight time.  }
         \label{fig:batch_traj_inference}
     \end{subfigure}
     \begin{subfigure}[th!]{\columnwidth}
         \centering
         \includegraphics[width=\columnwidth]{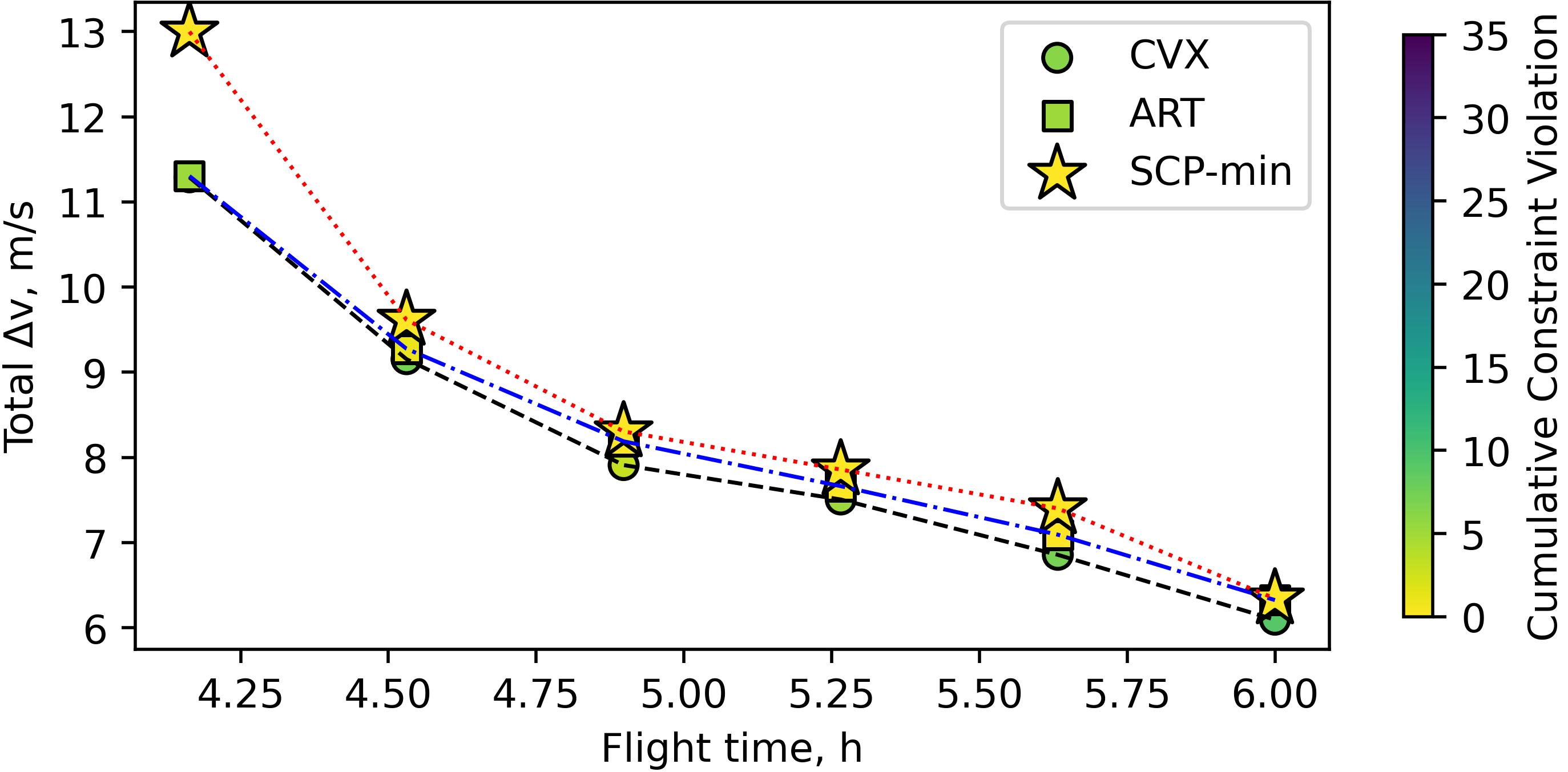 }
         \caption{Pareto-front of the competing objectives (flight time vs. fuel usage).}
         \label{fig:pareto_dv_tof}
    \end{subfigure} 
     \begin{subfigure}[th!]{\columnwidth}
         \centering
         \includegraphics[width=\columnwidth]{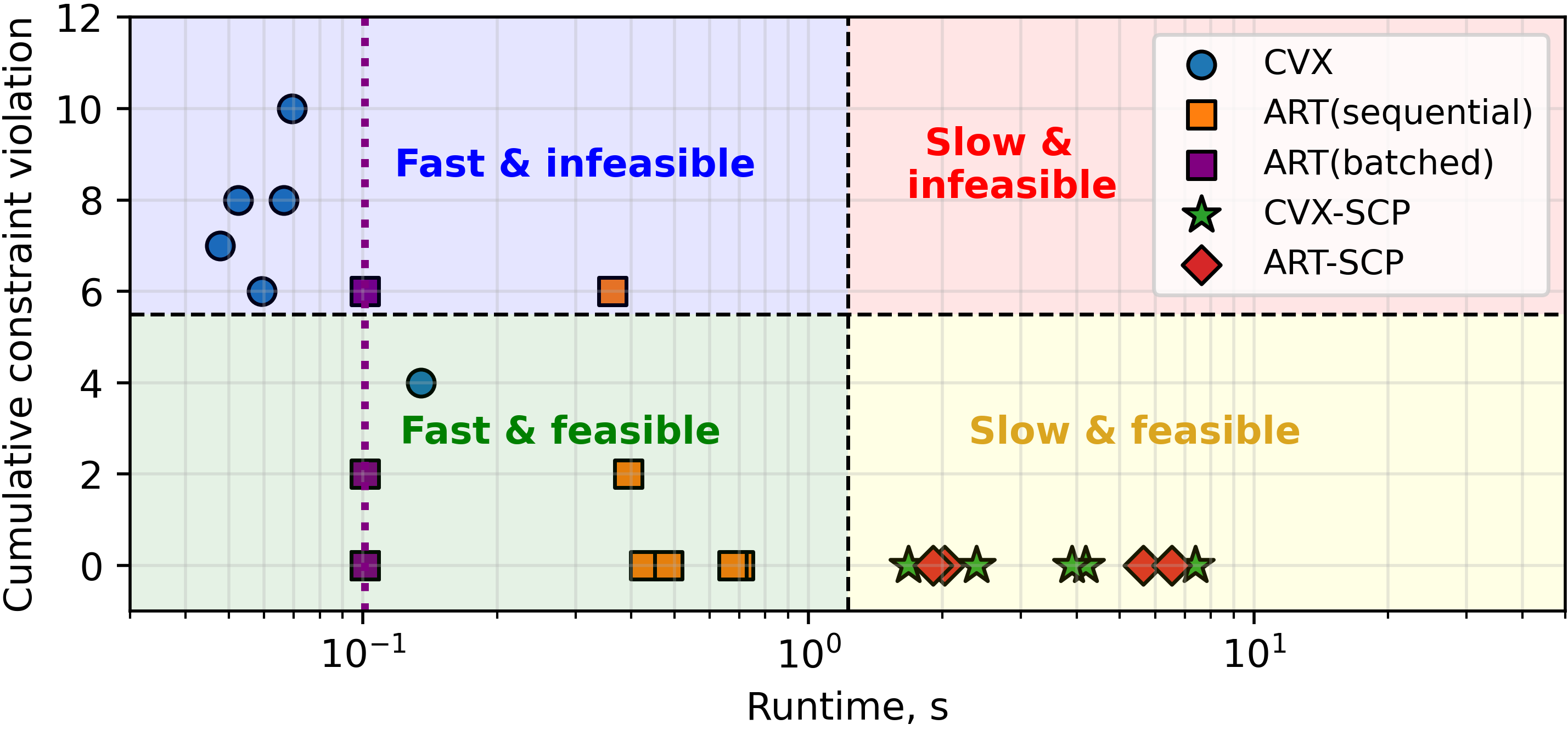}
         \caption{Distribution of constraint violation and the runtime of each trajectory generation method. 
         Each background color qualitatively indicates fast and feasible (green), slow and feasible (yellow), fast and infeasible (blue), and slow and infeasible (red). }
         \label{fig:pareto_ctg_runtime}
    \end{subfigure} 
    \caption{Near-Pareto-optimal trajectory generation via ART for a representative test scenario (cf. Eq.~\eqref{eq:rep_scenario}).}
    \label{fig:batch_inference}
\end{figure}

This subsection evaluates ART's effectiveness as an efficient and accurate surrogate for the Pareto front, as well as its practical utility within the proposed framework. 
Representative outputs from ART are shown in Figure~\ref{fig:batch_inference}.
The input consists of three components: (i) target orbital elements, (ii) boundary conditions (initial and terminal), and (iii) the range of admissible flight times over which the tradespace is explored (cf. Figure~\ref{fig:art_md_summary}). 
In particular, the figure summarizes the outputs of the following test scenario inputs, which are provided to a model trained for the far-range elliptic-orbit rendezvous case:
\begin{subequations} \label{eq:rep_scenario}   
\begin{align} 
\boldsymbol{\alpha}(t_0) & = [a \ [\text{km}],e,i\ [\text{rad}],\Omega \ [\text{rad}],\omega \ [\text{rad}],M \ [\text{rad}]] \nonumber \\
& = [9.50\cdot10^3, 0.195, 0.629, 0.0, 4.40, 2.51] , \\
a\boldsymbol{x}_{\text{e}}(t_0) & = a[\delta a, \delta\lambda, \delta e_x, \delta e_y, \delta i_x,\delta i_y] \ \text{[km]} \nonumber \\
& = [-15.1, -300,  53.1, -19.6, 3.05, -0.18] , \\
4.16  & \leq t_f \leq 6.0 \ \text{[h]}.
\end{align}
\end{subequations}
Based on this, the pipeline autoregressively generates a batch (here, $B=6$) of near-optimal trajectories spanning different flight times (cf. Figure~\ref{fig:batch_traj_inference}) together with their corresponding tradespace in the multi-objective domain (cf.  Figure~\ref{fig:pareto_dv_tof}). 
In Figure~\ref{fig:pareto_dv_tof}, ``SCP-min'' denotes the minimum cost obtained across the converged solutions of CVX-SCP and ART-SCP, serving as an approximation of the closest achievable true Pareto set.
These results demonstrate that ART consistently predicts the SCP cost more accurately than CVX.
This predictive accuracy yields an offset that allows for a conservative approximation of the Pareto surface, inherently accounting for future constraint satisfaction—a feature that relaxed optimization formulations cannot capture by design.

Furthermore, Figure~\ref{fig:pareto_ctg_runtime}  compares the runtime of trajectory generation (via inference or optimization) with the corresponding cumulative constraint violations.
The reported runtimes include only the optimization or inference stages, excluding preparatory steps such as the computation of STMs. 
Results are shown for both sequential and batched ART inference. To ensure a fair comparison, the batched inference time is normalized by the number of trajectories generated, since all trajectories within a batch are identical in terms of constraint violation.
The results indicate that even in the sequential setting, ART produces trajectories 3–10 times faster than SCP. 
When inference is parallelized, the computational cost of a single forward pass remains nearly constant with respect to batch size. 
As a result, the average per-trajectory runtime of batched inference approaches that of convex optimization, underscoring the efficiency and scalability of ART as a practical alternative to traditional optimization-based trajectory generation.

\begin{figure*}[ht!]
    \centering
     \begin{subfigure}{0.33\textwidth}
         \centering
         \includegraphics[width=\linewidth]{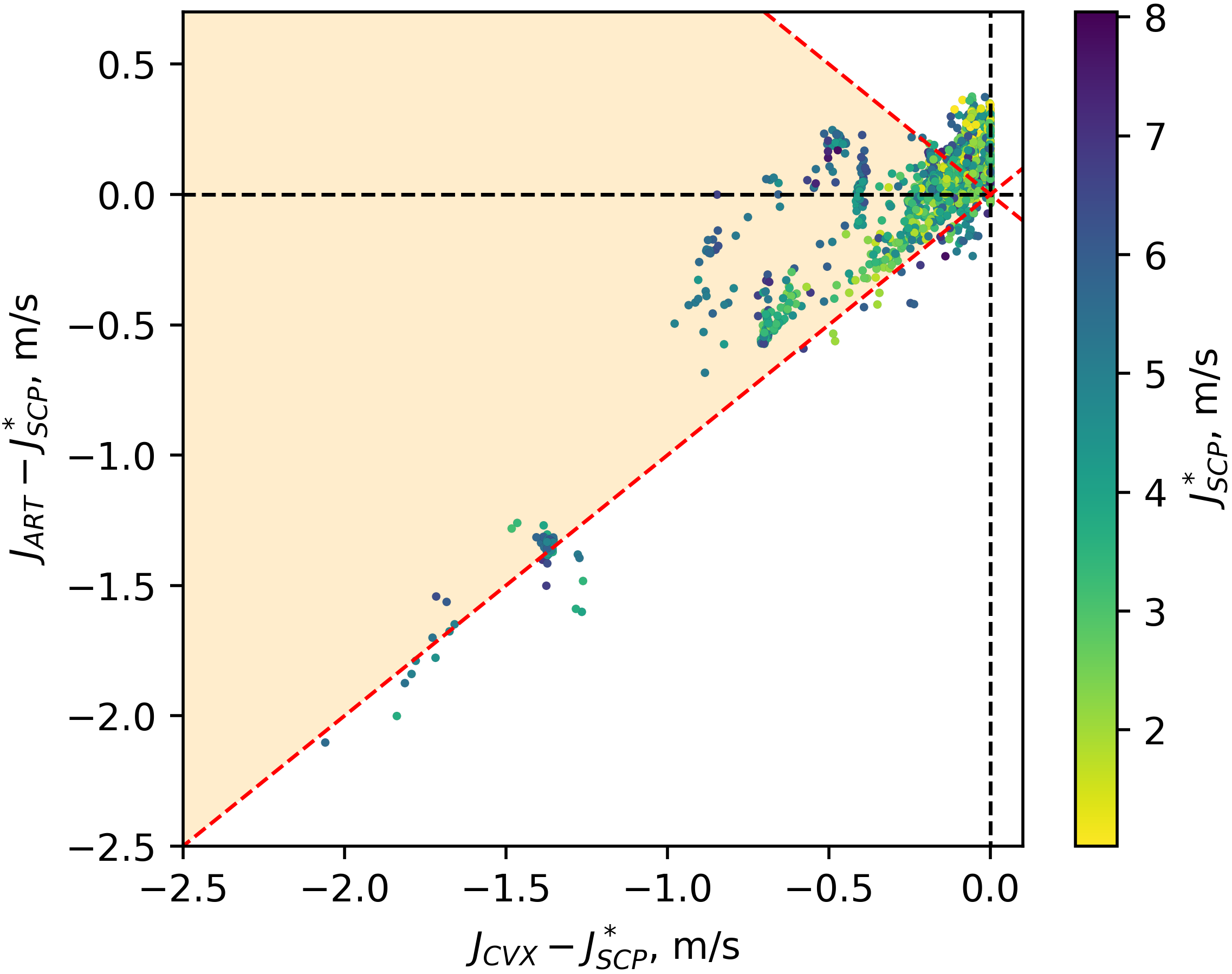}
         \caption{Generalization across time for the LEO scenario.}
         \label{fig:cost_offset_leo}
     \end{subfigure}
     \begin{subfigure}{0.33\textwidth}
         \centering
         \includegraphics[width=\linewidth]{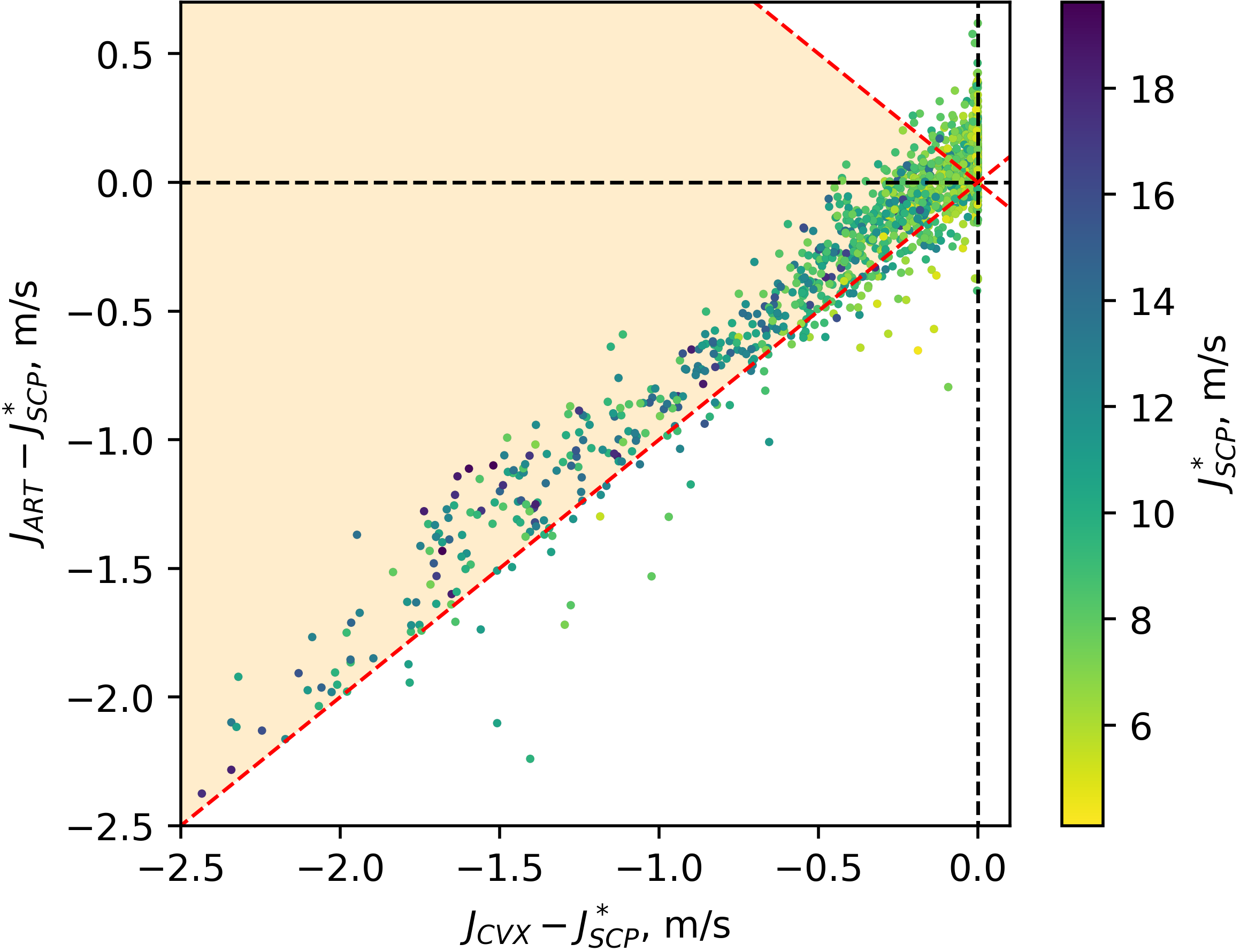}
         \caption{Generalization across time and dynamics for the elliptic orbit scenario.}
         \label{fig:cost_offset_kep_test}
    \end{subfigure} 
    \begin{subfigure}{0.33\textwidth}
         \centering
         \includegraphics[width=\linewidth]{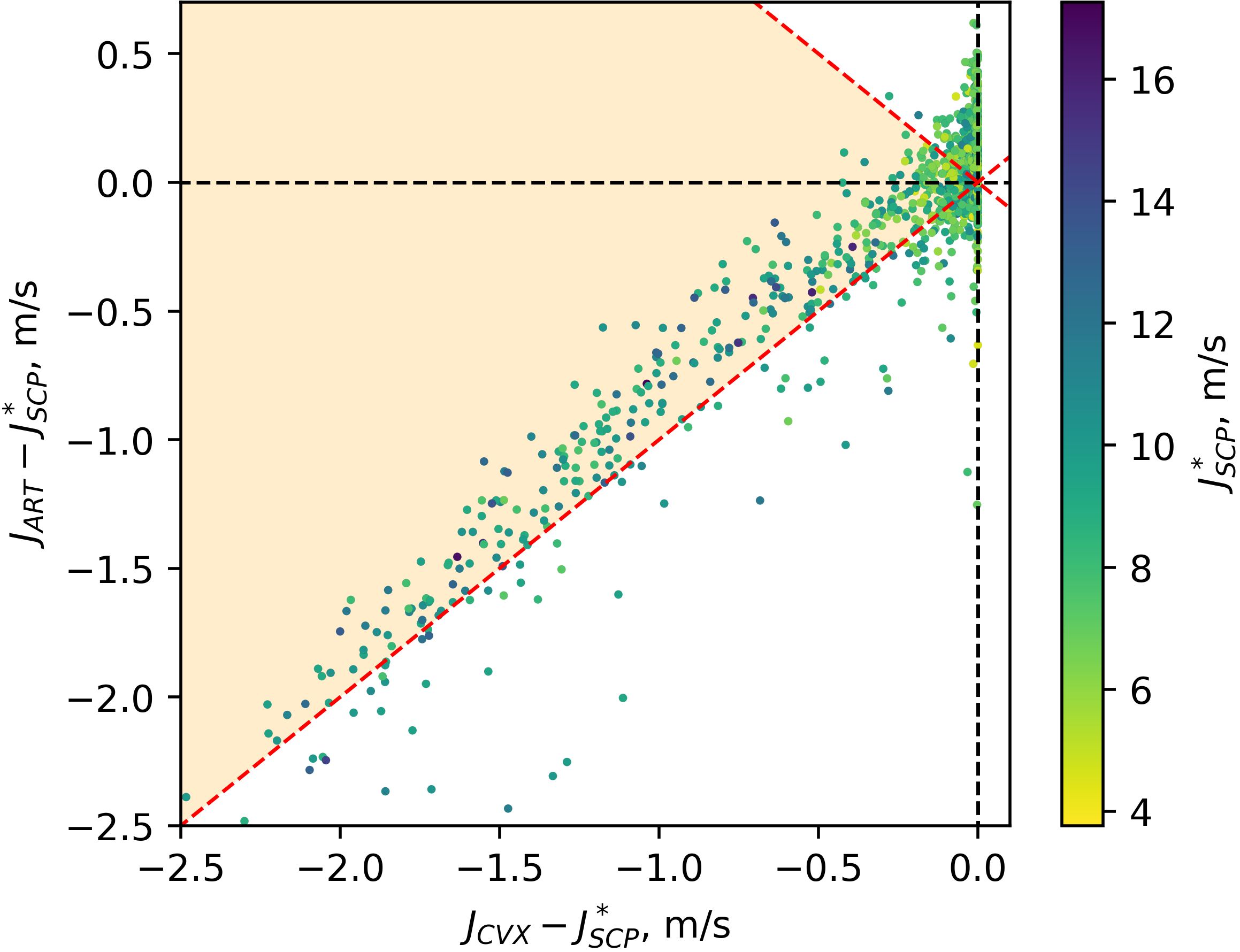}
         \caption{Generalization across time and dynamics with the elliptic orbit scenario according to the SpaceTrack debris dataset.}
         \label{fig:cost_offset_kep_debris}
    \end{subfigure} 
    \caption{Comparison of ART-inferred control costs and the CVX lower bound as mission-cost predictors (for a chosen flight time), shown relative to the minimum converged SCP objective. The orange region highlights cases where ART provides a more accurate cost prediction than CVX.}
    \label{fig:cost_offset}
\end{figure*}

The cost-prediction capability of ART is further evaluated through statistical analysis. 
Figure~\ref{fig:cost_offset} compares ART’s predicted control costs against those of CVX, using the performance analysis dataset introduced in the previous subsection (cf. Figure~\ref{fig:ws-analysis}). 
Both methods are evaluated relative to the minimum converged objective between CVX-SCP and ART-SCP.
Overall, ART demonstrates a consistent trend toward more accurate cost prediction across all three test cases, corroborating the observation from Figure~\ref{fig:pareto_dv_tof}. 
In the LEO scenario (cf. Figure~\ref{fig:cost_offset_leo}), when the offset of SCP from CVX is small, ART predictions remain close to the reference line $\mathcal{J}_{\text{ART}} = \mathcal{J}_{\text{SCP}}^* $, indicating that ART effectively captures how the CVX objective shifts once future constraint satisfaction is enforced. 
However, as the offset of SCP from CVX increases (i.e., $\mathcal{J}_{\text{SCP}} - \mathcal{J}_{\text{CVX}}^*$ grows), ART predictions exhibit a bias toward the CVX solution, with distributions clustering around $\mathcal{J}_{\text{ART}} = \mathcal{J}_{\text{CVX}}$. 
This behavior stems from the fact that the initial RTG is defined relative to the CVX objective.

For the elliptic orbit scenarios (Figs.~\ref{fig:cost_offset_kep_test} and \ref{fig:cost_offset_kep_debris}), the results show a nearly linear trend with an average offset of approximately 0.25 m/s in the direction of more accurate prediction relative to SCP. 
This improvement is particularly evident in scenarios where the gap between $\mathcal{J}_{\text{SCP}}^*$ and $\mathcal{J}_{\text{CVX}}^*$ is large, enabling ART to provide more reliable cost estimates of the SCP solution.
Although a slightly greater number of outliers are observed---cases in which ART predicts costs less accurately than CVX---the average cost offset in the test case with the debris data (cf. Figure~\ref{fig:cost_offset_kep_debris}) remains consistent with that obtained from the in-distribution evaluation (cf. Figure~\ref{fig:cost_offset_kep_test}).
This consistency highlights ART’s strong generalization capability. 
Overall, ART consistently improves the accuracy of cost estimation for the nonconvex OCP. 
Nevertheless, further refinement, such as direct prediction from mission design inputs, could enhance its predictive power and tighten the distribution around the reference line $\mathcal{J}_{\text{ART}} = \mathcal{J}_{\text{SCP}}^* $.

\section{Conclusion} \label{sec:conclusion}

This paper presents a rapid mission design framework supported by the multi-objective spacecraft trajectory generation under complex nonconvex constraints. 
A key enabler is the transformer model that is trained to predict trajectories conditioned on flight time and target orbital elements, allowing it to generate diverse, high-quality solutions across various mission scenarios. 
The trained model is validated against real orbital debris data spanning a wide range of orbital elements, demonstrating strong performance as an initial guess for chance-constrained, passively safe rendezvous trajectory optimization, suitable for subsequent refinement. 
Moreover, results show that the generated trajectories themselves can act as an effective surrogate model of the mission profile, particularly in challenging cases where convex relaxations severely violate nonconvex constraints.
By exploiting batched trajectory generation, the process achieves speeds approaching those of solving the convex relaxation itself.
Overall, the findings highlight the potential of AI-based trajectory generation as a powerful tool for agile trade-space exploration in preliminary mission design. 
The proposed framework is readily extensible to a broad range of applications, from six-degree-of-freedom rendezvous to interplanetary transfers. 
Future research will focus on mitigating residual terminal constraint violations arising from autoregressive generation, for example by incorporating multiple-shooting-based trajectory inference. 
In addition, further reductions in inference runtime would enable onboard mission design and trade-space analysis, potentially allowing these capabilities to become an integral component of future spacecraft autonomy.

\section*{Appendix} \label{appendix} 

\subsection{Transformer Model}
The presented transformer-based trajectory generation is implemented in PyTorch \cite{paszke2017automatic} and builds off Hugging Face's \texttt{transformers} library \cite{HuggFaceTransf}.
Table \ref{tab_appendix:hyper} presents the hyperparameter settings used in this work.

\begin{table}[ht]
\centering
\caption{Hyperparameters of ART architecture.}
\begingroup
\renewcommand*{\arraystretch}{1.25}
\begin{tabular}{l l}
    \thickhline
     Hyperparameter & Value \\
    \hline
     Number of layers & 6\\
     Number of attention heads & 6 \\
    Embedding dimension & 384 \\
     Batch size& 4 \\
    Context length $K$ & 50 \\
    Non-linearity & ReLU\\
    Dropout & 0.1\\
    Learning rate & 3e-5\\
    Grad norm clip & 1.0 \\
    Learning rate decay & None \\
    Gradient accumulation iters & 8 \\
    \thickhline
    \end{tabular}%
    \label{tab_appendix:hyper}
    \endgroup
\end{table}

\vspace{1em}

\subsection{SCP Parameters} \label{sec:appendix_scp}

Table \ref{tab_appendix:scp_param} presents the hyperparameters of the SCP used in the paper. 
The nomenclature corresponds to Oguri \cite{oguri2023successive}.

\begin{table}[ht]
\centering
\caption{Hyperparameters of $\texttt{SCVx}^*$.}
\renewcommand*{\arraystretch}{1.25}
\begin{tabular}{cc}
    \thickhline
    Hyperparameter & Value \\
    \hline
     $\epsilon$ & $10^{-3}$    \\
     $\{\rho_0, \rho_1, \rho_2\}$ & $\{ 0.0, 0.25, 0.7 \}$ \\
     $\{\alpha_1, \alpha_2, \beta, \gamma \}$ & $\{ 2,2,1.5,0.9\}$ \\
     $\{ r^{(1)}, r_{\min}, r_{\max} \}$ & $\{ 0.5, 10^{-6}, 10 \}$ \\
     $\{w^{(1)}, w_{\max}\}$ & $\{ 10, 10^9 \}$ \\
     $\#$ max. iter & 50  \\
    \thickhline
    \end{tabular}
    \label{tab_appendix:scp_param}
\end{table}

\vspace{1em}

\subsection{Uncertainty Modeling of Trajectory Optimization}

Three uncertainty sources are considered in the trajectory optimization: navigation error, control execution error, and unmodeled acceleration. 
The parameters for the uncertainty model are summarized in Table \ref{tab_appendix:uncertainty}.

\begin{table}[ht]
\centering
\caption{Uncertainty models}
\renewcommand*{\arraystretch}{1.25}\footnotesize
\begin{tabular}{cc}
    \thickhline
    Parameter & Value \\
    \hline
    Navigation error &   \\
    \hline
    $\boldsymbol{s}$ & 
    $\footnotesize{\begin{cases}
        \boldsymbol{s}_{1}
        & \text{if } \rho > \rho_1,\\
        \dfrac{\boldsymbol{s}_{1} (\rho - \rho_2) + \boldsymbol{s}_{2} (\rho_1 - \rho)}{\rho_1-\rho_2} 
        & \text{if } \rho_2 < \rho < \rho_1, \\
        \boldsymbol{s}_{2} 
        & \text{otherwise},
    \end{cases}}$ \\ 
    $\boldsymbol{s}_{1}$ &  $[4,400,4,2,2,4] \cdot 10^{-5} [\text{m}]$ \\
    $\boldsymbol{s}_{2}$ &  $[1,40,20,20,20,20] \cdot 10^{-4} [\text{m}]$ \\
    $\{\rho_1, \rho_2\} $ & $\{10^4, 10^3\} \ [\text{m}]$   \\
    \hline
    Control error \cite{gates1963simplified} & \\
    \hline 
    $\{\sigma_{m1}, \sigma_{m2}\}$  &  $\{2 \ [\text{mm/s}], 3\cdot10^{-4}\}$   \\
    $\{\sigma_{p1}, \sigma_{p2}\}$  &  $\{0.3 \ [\text{mm/s}], 3\cdot10^{-4}\}$   \\
     \hline 
    Unmodeled acc. & \\
    \hline 
    $Q_{ki}$ &  $ 10^{-6}\cdot \boldsymbol{I}_6 \ [\text{m}^2]$\\
    \hline 
    Risk factor & \\
    \hline 
    $\Delta_k$ & 0.003 ($3\sigma$-safety) \\
    \thickhline
    \end{tabular}
    \label{tab_appendix:uncertainty}
\end{table}

\section*{Acknowledgements}

This work is supported by Blue Origin (SPO \#299266) as an Associate Member and Co-Founder of Stanford’s Center of AEroSpace Autonomy Research (CAESAR). 
This article solely reflects the opinions and conclusions of its authors and not any Blue Origin entity. 
The authors are grateful to Davide Celestini for his valuable insights and to Sukeerth Ramkumar for his assistance.
Yuji Takubo acknowledges financial support from the Ezoe Memorial Recruit Foundation.

\bibliographystyle{IEEEtran}
\bibliography{reference}

\input{bio}

\end{document}

%% file: bio.tex
\thebiography

\begin{biographywithpic}
{Yuji Takubo}{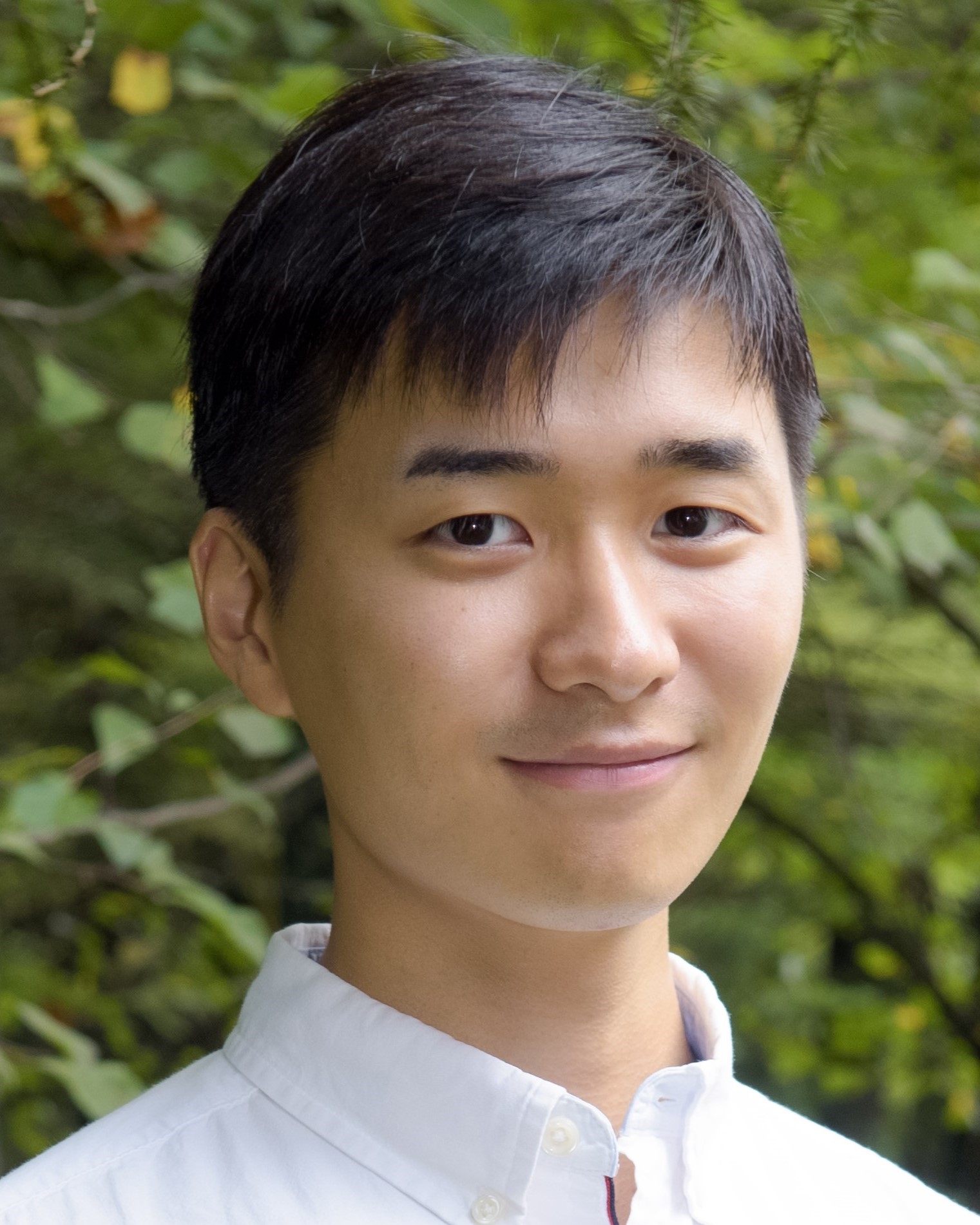}
Yuji Takubo is a Ph.D. candidate in the Space Rendezvous Lab at Stanford University. He received the M.S. in Aeronautics and Astronautics from Stanford University (2025) and the B.S. in Aerospace Engineering from the Georgia Institute of Technology (2023). 
He is a recipient of the NASA Jet Propulsion Laboratory Graduate Fellowship (2023) and the Ezoe Memorial Recruit Foundation Fellowship (2020-).  
His research focuses on optimization and astrodynamics, particularly with emphasis on distributed space systems. 
\end{biographywithpic}

\begin{biographywithpic}
{Daniele Gammelli}{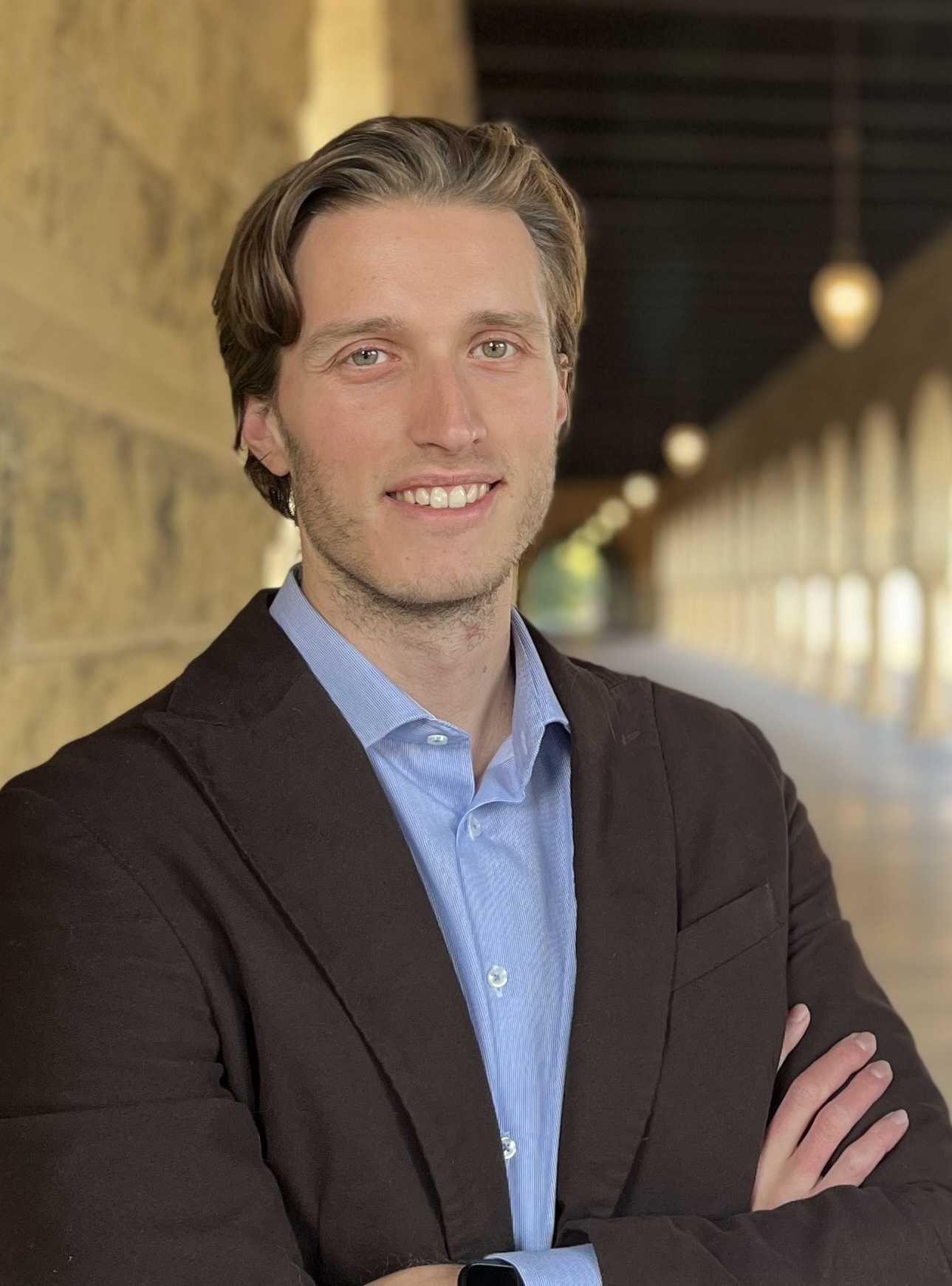}
Dr. Daniele Gammelli is a Postdoctoral Scholar in the Autonomous Systems Lab at Stanford University. He received the Ph.D. in Machine Learning and Mathematical Optimization at the Department of Technology, Management and Economics at the Technical University of Denmark. Dr. Gammelli’s research focuses on developing learning-based solutions that enable the deployment of future autonomous systems in complex environments, with an emphasis on large-scale robotic networks, aerospace systems, and future mobility systems. During his doctorate and postdoctorate career, Dr. Gammelli has been making research contributions in fundamental AI research, robotics, and its applications to network optimization and mobility systems. His research interests include deep reinforcement learning, generative models, graph neural networks, bayesian statistics, and control techniques leveraging these tools.
\end{biographywithpic}

\begin{biographywithpic}
{Marco Pavone}{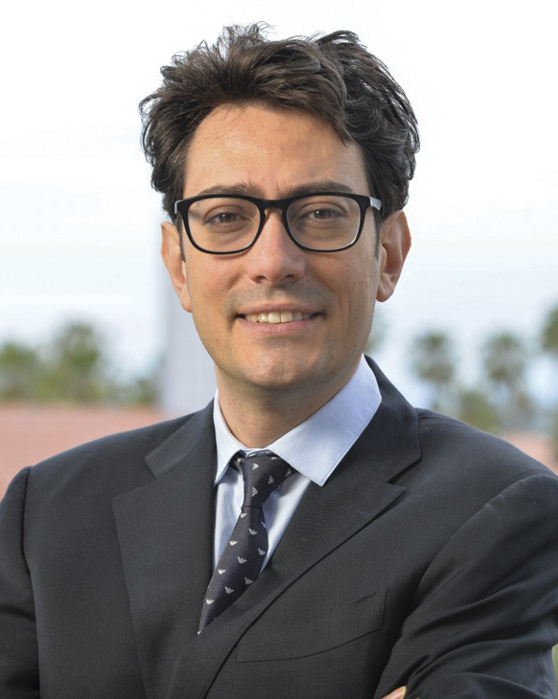}
Dr. Marco Pavone is an Associate Professor of Aeronautics and Astronautics at Stanford University, where he directs the Autonomous Systems Laboratory and the Center for Automotive Research at Stanford. He is also a Distinguished Research Scientist at NVIDIA where he leads autonomous vehicle research. Before joining Stanford, he was a Research Technologist within the Robotics Section at the NASA Jet Propulsion Laboratory. He received a Ph.D. degree in Aeronautics and Astronautics from the Massachusetts Institute of Technology in 2010. His main research interests are in the development of methodologies for the analysis, design, and control of autonomous systems, with an emphasis on self-driving cars, autonomous aerospace vehicles, and future mobility systems. He is a recipient of a number of awards, including a Presidential Early Career Award for Scientists and Engineers from President Barack Obama, an Office of Naval Research Young Investigator Award, a National Science Foundation Early Career (CAREER) Award, a NASA Early Career Faculty Award, and an Early-Career Spotlight Award from the Robotics Science and Systems Foundation. He was identified by the American Society for Engineering Education (ASEE) as one of America's 20 most highly promising investigators under the age of 40. His work has been recognized with best paper nominations or awards at a number of venues, including the European Conference on Computer Vision, the IEEE International Conference on Robotics and Automation, the European Control Conference, the IEEE International Conference on Intelligent Transportation Systems, the Field and Service Robotics Conference, the Robotics: Science and Systems Conference, and the INFORMS Annual Meeting.
\end{biographywithpic}

\begin{biographywithpic}
{Simone D'Amico}{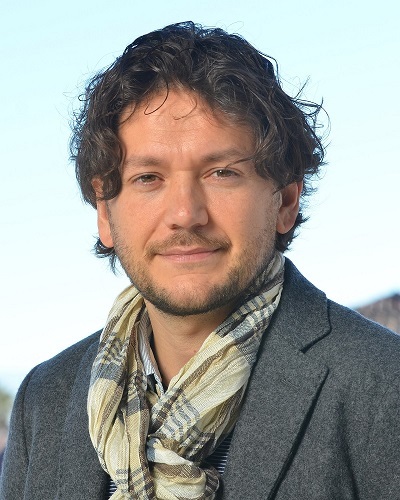}
Dr. Simone D'Amico is Associate Professor of Aeronautics and Astronautics (AA) and Professor of Geophysics (by Courtesy). He is the Founding Director of the Stanford Space Rendezvous Laboratory, Director of the Undergraduate Program in AA, and Founding Co-Director of the Center for AEroSpace Autonomy Research (CAESAR) at Stanford. He is a W.M. Keck Faculty Scholar in the School of Engineering and Hoover Science Fellow.
He has 20+ years of experience in research and development of autonomous spacecraft and distributed space systems, including multi-agent architectures such as “rendezvous, proximity operations, and capture”, “satellite formation-flying and swarms”, “fractionated spacecraft”, and “mega-constellations”. He developed the GNC system of several such missions (e.g., TanDEM-X, PRISMA) and is currently the institutional PI of four autonomous satellite swarms funded by NASA and NSF with one of them (NASA Starling) operational in orbit.
Besides academia, Dr. D’Amico is the Chief Science Officer (CSO) of EraDrive. He is in the Advisory Board of four space start-ups focusing on distributed space systems for future applications in SAR remote sensing, orbital lifetime prolongation, and space-based solar power. He received the B.S. and M.S. degrees from Politecnico di Milano (2003) and the Ph.D. degree from Delft University of Technology (2010). Before Stanford, Dr. D’Amico was research scientist and team leader at the German Aerospace Center (DLR). He was the recipient of several awards, including the 2020 IEEE M. Barry Carlton Award, Best Paper Awards at AMOS (2025), IAF (2022), IEEE (2021), AIAA (2021), AAS (2019) conferences, the Leonardo 500 Award by the Leonardo da Vinci Society/ISSNAF (2019), FAI/NAA’s Group Diploma of Honor (2018), DLR’s Sabbatical/Forschungssemester (2012) and Wissenschaft Preis (2006), and NASA’s Group Achievement Awards for the STARLING (2024) and GRACE mission (2004).
\end{biographywithpic}